\documentclass[12pt,leqno]{article}
\usepackage{color}
\usepackage{enumitem}

\def\int{\mathbb{Z}}

\def\gg{\mathfrak g}

\def\JJ{{\mathfrak{J}}}
\def\ll{{\mathfrak l}}
\def\hh{{\mathfrak h}}

\def\WW{{\mathcal W}}
\def\WK{\WW_{\HH,K}}
\def\WKl{\WW_{\HH,K}^{l}}
\def\WHl{\WW_{\HH}^l}
\def\WKl{\WW_{\HH,K}^l}

\def\cc{\mathfrak{c}}

\def\zz{\mathfrak{Z}}

\def\HH{\mathcal{H}}
\def\HHl{\mathcal{H}_l}

\def\AA{{\mathcal{A}}}
\def\Lc{L_{\mathbb C}}
\def\Ec{E_{\mathbb C}}
\def\Vc{V_{\mathbb C}}
\def\Abar{\overline{\AA}}
\def\Abarp{\overline{{\mathcal A}'}}
\def\Abarc{{\overline{\AA}}_{\mathbb C}}
\def\WH{\WW_\HH}
\def\hWH{\widehat{\WW}_\HH}
\def\HHc{\HH_{\mathbb C}}

\def\pf{\proof}
\def\id{{\rm id}}

\def\epf{\hfill$\Box$ \medskip}
\def\interiorc{\mathring{\Lc}}
\def\interior{\mathring{L}}
\def\interiorcp{\mathring{L'_{\mathbb C}}}
\def\interiorp{\mathring{L'}}
\usepackage{tikz}
\usepgflibrary{shapes.geometric}
\usepackage{amsmath,amsthm, amscd, amssymb, amsfonts, mathrsfs}
\usepackage[all]{xy}
\usepackage{color}
\usepackage[vcentermath]{youngtab}

\usepackage{graphics}
\usepackage{amssymb}
\usepackage{amscd}
\usepackage{amsmath}
\input xy 
\xyoption{all}
\title{Affine hyperplane arrangements\\ and Jordan classes}
\newtheorem{theorem}{Theorem}[section]
\newtheorem{lemma}[theorem]{Lemma}
\newtheorem{corollary}[theorem]{Corollary}
\newtheorem{proposition}[theorem]{Proposition}
\newtheorem{definition}[theorem]{Definition}
\newtheorem{remark}[theorem]{Remark}
\newtheorem{example}[theorem]{Example}

\author{Giovanna Carnovale, Francesco Esposito\\
Dipartimento di Matematica ``Tullio Levi-Civita''\\
Torre Archimede - via Trieste 63 - 35121 Padova - Italy\\
email: carnoval@math.unipd.it, esposito@math.unipd.it }
\date{}

\begin{document}
\maketitle
\begin{abstract}
We study the geometry of the stratification induced by an affine hyperplane arrangement $\HH$ on the quotient of a complex affine space by the action of a discrete group  preserving $\HH$. We give conditions ensuring normality or normality in codimension $1$ of strata. 
% In special cases this stratification is equivalent to the stratification of categorical quotients of closures of Jordan classes in a complex semisimple algebraic group or  Lie algebra. 
As an application, we provide the list of the categorical quotients of closures of Jordan classes and of sheets in all complex simple algebraic groups that are normal. In the simply connected case, we show that normality of such a quotient is equivalent to its smoothness. 
\end{abstract}


\section{Introduction}

In \cite{BK,bo} the stratification of a semisimple Lie algebra by Jordan classes (also called decomposition classes or packets) was introduced and studied in order to describe the sheets for the adjoint action of a semisimple algebraic group $G$  on its Lie algebra $\gg$. It was shown that every sheet is the regular part in the closure of a unique Jordan class and as a consequence, sheets could be classified in terms of combinatorial data. 
Closure relations for Jordan classes were explicitly given and, for $\mathcal S$ a sheet, the  topology of the orbit space ${\mathcal S}/G$  and 
%criteria were given in order to ensure that the $G$-module structure of ${\mathbb C}[\Oc]$, for adjoint orbits $\Oc$, is constant along sheets. 
the normalisation of the categorical quotient $\overline{\mathcal S}/\!/G$ were explicitly described. These quotients are the closures of Luna strata for $\gg/\!/G$, as defined in \cite[III.2]{luna}.

Subsequently, it was proved in  \cite{katsylo} that the orbit space  ${\mathcal S}/G$ can be given the structure of a  geometric quotient which  is isomorphic to the quotient of an affine space modulo the action of a finite group. Richardson in \cite{richardson} has provided a criterion ensuring normality of $\overline{\mathcal S}/\!/G$, where ${\mathcal S}$ is a sheet or a Jordan class and has produced a complete list of the normal quotients for classical Lie algebras. The list for exceptional Lie algebras was obtained in \cite{broer} and \cite{DR} with different techniques.
% The family of quotients $\overline{S}/\!/G$ where $S$ runs through the set of sheets in $\gg$ is the same as the one in which $S$ is a Jordan class containing semisimple elements. 
The same approach allowed to provide in \cite{richardson,broer,DR} the complete list of those regular Jordan classes whose closure is normal and Cohen-Macaulay.

In the seminal paper \cite{lusztig-inventiones}, Lusztig introduced a stratification on $G$ which is analogous to the partition into Jordan classes and proved that topological properties of this stratification (and of its quotient) encoded representation theoretic information for $G$ and its Weyl group. This stratification in Jordan classes for $G$ has been crucial for the study of sheets for the adjoint action of $G$  on itself \cite{gio-espo}. An analogue of Katsylo's result \cite{katsylo} for ${\mathcal S}/G$, where ${\mathcal S}$ is a sheet in $G$ consisting of spherical conjugacy classes was given in \cite{gio-espo-RT}. Other properties of quotients of sheets and of closures of Jordan classes, including the description of the normalisation of  
$\overline{\mathcal S}/\!/G$ in the spirit of \cite{bo} were given in \cite{proceedings-joseph}. These quotients are the irreducible components of closures of Luna strata for $G/\!/G$.
 The origin and main motivation of this paper is to detect when they are normal.
%normality of the quotients $\overline{J}/\!/G$, where $J$ is a Jordan class or a sheet in $G$. %As in the Lie algebra case, this family of quotients coincides with the family of quotients $\overline{S}/\!/G$, where $S$ is a sheet in $G$. 

In the Lie algebra case, Douglass and R\"ohrle \cite{DR} translated the normality condition in \cite{richardson} in terms of properties of subarrangements of the Weyl group hyperplane arrangement. We bring affine hyperplane arrangements into the picture in the following way. 
For a semisimple group $G$ we consider the stratification on a Cartan subalgebra $\hh$ in ${\rm Lie}(G)$ induced by the the corresponding affine Weyl group arrangement, and the quotient $\hh/\WW$ of $\hh$ by the action of a finite extension $\WW$ of the affine Weyl group. A Jordan stratum in $G/\!/G$  is shown to be (analytically) isomorphic to a stratum for the quotient stratification on $\hh/\WW$.  
In order to provide a uniform treatment of the Lie algebra case, the simply connected group case and the non simply connected group case, we put this problem into a more general framework, allowing a wider range of choices for the acting group $\WW$, which will be the extension of a group $\WH$ generated by suitable affine reflections by a finite group $K$. The quotients studied in \cite{DR} correspond to the case in which $\WH$ is finite and $K$ is trivial. Similar questions have been addressed in \cite{DR2}, by considering the action of finite complex reflections groups (with no extensions). 

One of the main novelties in our approach is the analysis of these stratifications by looking at the local geometry of strata. Around unibranch points strata are smoothly equivalent to strata for a finite arrangement. We describe how the geometry of the problem behaves along strata and reduce the verification of normality of a stratum $X$ to checking normality at  a well-chosen point in each minimal stratum contained in $X$. In the case of Jordan classes in semisimple groups, minimal strata correspond to the Jordan classes consisting of one single  conjugacy class, which is necessarily isolated, in the terminology of \cite{lusztig-inventiones}. Around such points the stratum is smoothly equivalent to the quotient of the closure of a Jordan class in a Lie algebra with automorphisms.

 %By a result of Richardson,  $X_J$ is normal if and only if the closure of the regular Jordan class $J'$ with semisimple part $J$ is normal, so our list provides also a list of regular Jordan classes that are normal. DA FARE ANCORA. 

Our approach sheds light on some phenomena which could be observed in the previously cited papers. Most evidently, it explains in terms of normality in codimension one a rigidity property of the combinatorial data associated with normal quotients that was given in terms of equality of two families of exponents in \cite{DR,DR2}.  Our interpretation is obtained by associating to each quotient a $K$-stable family of faces $\Sigma(\HH,K,L)$ of a fundamental domain for the action of $\WH$. This set is a combinatorial counterpart for some geometric properties, e.g., both unibranchedness and normality can be read-off from the  properties of $\Sigma(\HH,K,L)$.   

As a final output we produce the list of normal strata for any simple $G$. We prove that when $G$ is simply connected, a stratum is normal if and only if it is smooth. The same phenomenon occurs in the Lie algebra case \cite{broer}. For finite complex reflection groups this was observed in \cite{DR2}. The results obtained here have been applied in \cite{ACE} to produce the complete list of regular Jordan classes in $G$ semisimple and simply connected whose closure is  normal and Cohen-Macaulay.

\section{Notation, basic definitions and motivation}
\subsection{Index of notation}
$\HH$ (admissible) affine hyperplane arrangement in $E$, page \pageref{HH}\\
$\WH$ group generated by the reflections with respect to affine  hyperplanes in $\HH$, page \pageref{HH}\\
$V(H)$, the direction of an affine hyperplane, page \pageref{HH}\\
$V_{\HH}=\bigcap_{H\in\HH}V(H)$, page \pageref{HH}\\
$E^{ess}$, the affine space on which $\WH$ is essential, page \pageref{ess}\\
$\HH^{ess}$, the restriction of $\HH$ to $E^{ess}$, page \pageref{ess}\\ 
$\Abar$ fundamental domain for the $\WH$-action on $E$, page \pageref{AA}\\
$K$ subgroup of ${\rm Stab}_{{\rm Aut}(E)}(\Abar)$, page \pageref{WK}\\
$\WK=K\ltimes\WH$, page \pageref{WK}\\
$\HHl=\{H\in\HH~:~l\in H\}$, page \pageref{Hl}\\
%$\WHl={\rm Stab}_{\WH}(l)$, page \pageref{Hl}\\
$L$ a flat, i.e., an intersections of affine hyperplanes in $\HH$, page \pageref{flat}\\
$\mathring{L}=L\setminus\bigcup_{L'\subsetneq L}L'$, page \pageref{flat}\\
$\HH^L=\{H\cap L~:~\emptyset\neq H\cap L\neq L\}$ induced hyperplane arrangement on $L$, page \pageref{flat}\\
${\mathcal C}(\HH)$ the set of chambers induced by $\HH$, page \pageref{CH}\\
${\mathcal P}(\HH)$ the poset of chambers, page \pageref{CH}\\
$F$ a face in ${\mathcal P}(\HH)$, page \pageref{CH}\\
$|F|$ affine subspace generated by $F$, page \pageref{CH}\\
$\WW^p$ stabiliser of the point $p$ in the group $\WW$, page \pageref{Wp}\\
$A_{\mathbb C}$, complexification of the affine or vector space $A$, page \pageref{complex}\\
$\Abarc$, fundamental region for the $\WH$-action on $E_{\mathbb C}$, page \pageref{complex-fundamental}\\ 
$X(\HH,K,L):=\WK\Lc/\WK$, a stratum in $\Ec/\WK$, page \pageref{XL}\\
$X(\HH,L):=\WH\Lc/\WH$, a stratum in $\Ec/\WH$, page \pageref{XL}\\
$e=\exp(2\pi i-)\colon \gg\to G_{sc}$, $e_\pi=\exp(2\pi i-)\colon \gg\to G$, page \pageref{exp}\\
$\sim_{se}$, smooth equivalence, page \pageref{se}\\
$\Gamma_{\HH,L}={\rm Stab}_{\WH}(L)$, page \pageref{gamma}\\
$\Gamma_{\HH,K,L}={\rm Stab}_{\WK}(L)$, page \pageref{gamma}\\
$\iota_L\colon\Gamma_{\HH,K,L}\to{\rm Aut}(L)$, the map induced by restriction to $L$, page \pageref{iotaL}\\
$\Gamma_{W,L}={\rm Stab}_{\WH}(e(\Lc))$, page \pageref{gamma}\\
$\tilde{X}(\HH,L)=\Lc/\Gamma_{\HH,L}$, page \pageref{gamma}\\
$\tilde{X}(\HH,K,L)=\Lc/\Gamma_{\HH,K,L}$, page \pageref{gamma}\\
$\Omega(L,\HH,K)_l:=\left\{w\Lc~:~ w\in\WK,\; l\in w\Lc\right\}$, page \pageref{omega}\\
$\WW_{\HH}^{L}:=\cap_{l\in\Lc}\WW_\HH^l$, page \pageref{WHKL}\\
$\WW_{\HH,K}^{L}:=\cap_{l\in\Lc}\WKl$, page \pageref{WHKL}\\
$K^{L}:=\cap_{l\in\Lc}{K}^{l}$, page \pageref{WHKL}\\
$U_L:=\{l\in\Lc~:~\WKl=\WW_{\HH,K}^{L}\}$, equation \eqref{eq:ul}\\
$Y(\HH,K,L)=\WK \Lc$, page \pageref{Y}\\
$Y(\HH,K,L)_l=\bigcup_{\genfrac{}{}{0pt}{}{w\in \WK,} 
{l\in w\Lc}}w\Lc$, page \pageref{Y}\\
$\Sigma(\HH,K,L):=\{w\interior\cap\Abar~:~w\in \WW_{\HH,K},\;\;wL \mbox{ lies over }\Abar \}$, page \pageref{Sigma}\\
  $\Sigma(\HH,L):=\{w\interior\cap\Abar~:~w\in \WW_{\HH},\;\;wL \mbox{ lies over }\Abar \}$, page \pageref{Sigma}\\
  ${\mathcal P}(\Sigma(\HH,K,L)):=\{F\in{\mathcal P}(\HH^L)~:~F\subset\overline{F'}, \text{ for some } F'\in\Sigma(\HH,K,L) \}$, page \pageref{posets}\\
${\mathcal P}(\Sigma(\HH,L)):=\{F\in{\mathcal P}(\HH^L)~:~F\subset\overline{F'} \text{ for some } F'\in\Sigma(\HH,L) \}$, page \pageref{posets}\\
${\mathcal F}_l=\{F\in\Sigma(\HH,K,L)~:~l\in\overline{F}\}$, page \pageref{Fl}\\
 $S$, the set of  Coxeter generators of $\WH$, page \pageref{S}\\
 $\{N_1,\,\ldots,\,N_n\}$, nodes of the Coxeter graph of $\WH$, page \pageref{S}\\
 $\{x_0,\,\ldots,\,x_n\}$, vertices of a fundamental alcove, page \pageref{vertices}\\
 $\omega_i^\vee$, $i=1,\,\ldots,\,n$, fundamental co-weights, page \pageref{coweights}\\
 $M_F$, the walls of $\Abar$ containing the face $F$, page \pageref{emmeffe}\\
 $S_F$, the subset of $S$ containing the reflections with respect to $H\in M_F$, page  \pageref{emmeffe}\\
 $\hat{S}_F=S\setminus S_F$, page  \pageref{emmeffe}\\
 %${\mathcal S}_{\HH,K,L}:=\{J\subset S~:~J=J_F \textrm{ for some }F\in\Sigma(\HH,K,L)\}$, page \pageref{SX}\\
 $[S']:=\{S''\subset S~:~wS'=S'' \textrm{ for some } w\in\WH\}$, the Coxeter class of $S'\subset S$, page \pageref{CoxeterJ}\\
 $I:=K/K'$ for $K'\lhd K$, page \pageref{inertia}\\
$I^{X'}={\rm Stab}_I(X(\HH, K',L))$, page \pageref{inertia}\\
$I^{X',x}:=I^x \cap I^{X'}$ for $x\in X(\HH,K',L)$, page \pageref{inertia}\\

\subsection{Hyperplane arrangements and main question}\label{sec:basic-hyperplanes}

\subsubsection{Basic definitions}\label{subsec:basic}Let $E$ be an Euclidean space with direction vector space $V$ acting simply transitively by translations on $E$. We denote by $\HH$  a (not necessarily finite) affine hyperplane arrangement in $E$ and by $\WH$  the group generated by the reflections with respect to the affine hyperplanes in $\HH$.\label{HH} We say that $\HH$ is admissible if $\WH$, equipped with the discrete topology, acts properly on $E$ and preserves $\HH$. In this case, $\HH$ is locally finite \cite[V.3.1]{bourbaki}. 

We denote by $V(L)\subset V$ the direction of an affine subspace $L\subset E$ and we set $V_{\HH}=\cap_{H\in\HH}V(H)\subset V$. The action of the group $\WH$ and $\HH$ are called essential if $V_{\HH}=\{0\}$. 
The inclusion $\WH\leq {\rm Aut}(E)\simeq V\rtimes {\rm GL}(V)$ followed by the natural quotient by $V$ gives a linear $\WH$-action on $V$. Although the isomorphism  ${\rm Aut}(E)\simeq V\rtimes {\rm GL}(V)$ depends on the choice of a point, different choices induce isomorphic $\WH$-actions, all fixing $V_{\HH}$ pointwise. 
%The natural projection $\pi\colon E\to E/V_{\HH}$ is $\WH$-equivariant
%and gives rise to a principal $V_{\HH}$-bundle.  
%and 
For any $p\in E$ the affine subspace $E^{ess}=p+V_{\HH}^\perp$ is $\WH$-stable because $\WH$ is generated by (affine) reflections in directions that are orthogonal to $V_{\HH}$ and we have a $\WH$-equivariant isomorphism $E^{ess}\simeq E/V_\HH$.
%and a $\WH$-equivariant decomposition $E=E^{ess}+V_{\HH}$.
%gives a $\WH$-equivariant trivialisation of $\pi$. 
In addition, if $\HH^{ess}$ is the arrangement induced on $E^{ess}$, then it is essential, $\WH\simeq\WW_{\HH^{ess}}$ and $\HH^{ess}$ is admissible if $\HH$ is so. \label{ess}

The intersections of hyperplanes in $\HH$ induce a stratification on $E$, whose parts are called flats.\label{flat} For a flat $L$, we set $\mathring{L}:=L\setminus\bigcup_{L'\subsetneq L}L'$ and denote by $\HH^L$ the hyperplane arrangement induced by $\HH$ on $L$. In general, $\HH^L$ is not admissible even if $\HH$ is so. 

The connected components of $E\setminus\bigcup_{H\in\HH}H$ are called the chambers of $\HH$: we denote by ${\mathcal C}(\HH)$ the set of chambers of $\HH$ and similarly ${\mathcal C}(\HH^L)$ the set of chambers of $\HH^L$ in $L$. The closure of a chamber is a convex polytope. Following \cite[Definition 2.18]{OT} we set \label{CH}
\begin{align*}{\mathcal P}(\HH)=\bigcup_{L\text{ a flat in $E$}}{\mathcal C}(\HH^L)\end{align*} and we view it as a collection of subsets of $E$. Any $F\in  {\mathcal P}(\HH)$ is called a face, $|F|$ will denote the support of $F$, i.e. the minimal affine space containing $F$. Each face is open in its support and we set $\dim F:=\dim |F|$. By $\overline{F}$ we usually mean the closure in $E$. The set ${\mathcal P}(\HH)$ has a natural poset structure given by inclusion of closures, we shall call it the face poset of $\HH$. For a given chamber $C$, we say that $H\in\HH$ is a wall of $C$ if $H$ is the support of a (maximal) face of $C$. For any subposet in ${\mathcal P}(\HH)$, a gallery is a sequence of equidimensional faces $F_i$ for $i=0,\,\ldots,\,\ell$ such that for every $i$ there is a unique face $F'_i$ of codimension $1$ contained in $\overline{F_i}\cap\overline{F_{i+1}}$. 

{\em From now on we assume that $\HH$ is admissible.} Under this assumption $\WH$ acts simply transitively on ${\mathcal C}(\HH)$, \cite[Th\'eor\`eme 1, V.3.2]{bourbaki} and the closure of a chamber is a fundamental domain for the action of $\WH$, \cite[Th\'eor\`eme 2, V.3.3]{bourbaki}. We fix such a fundamental domain $\overline{\AA}$. If $\HH$ is not essential we take $\overline{\AA}=\overline{\AA}^{ess}+V_{\HH}$, where $\overline{\AA}^{ess}$ is a fundamental domain for the action of $\WH$ on $E^{ess}=E/V_\HH$.\label{AA}

By \cite[V.3.7, V.3.8]{bourbaki}, if the $\WH$-action on $V$ is not irreducible, then  $V=V_{\HH}\oplus(\bigoplus_{j=1}^rV_{(j)})$, where each $V_{(j)}$ for $j>0$ is an irreducible representation of $\WH$, $\WH\simeq\prod_{i=j}^r{\WH}_{(j)}$ where each factor ${\WH}_{(j)}$ acts irreducibly on $V_{(j)}$ and trivially on  $V_{(j')}$ for $j'\neq j$ and is generated by the reflections with respect to the hyperplanes in the induced, admissible arrangement $\HH_{(j)}$ on $E_{(j)}=p+V_{(j)}$, for some $p\in E$. Hence, $E\simeq V_\HH\times(\prod_{j=1}^rE_{(j)})$. 
%
%
%then $\WH$, $V$ and $E$ decompose as $\WH=\prod_{i=j}^r{\WH}_{(j)}$,  and $E=V_0(\prod_{j=1}^rE_{(j)})$ so that  for $j>0$  each factor ${\WH}_{(j)}$ acts irreducibly on $V_{(j)}$ and trivially on  $V_{(j')}$ for $j'\neq j$, it is %generated by reflections with respect to the affine hyperplanes in the induced arrangement $\HH_{(j)}$ on $E_{(j)}=p+V_{(j)}$, and ${\WH}_{(j)}$-action is admissible, irreducible and effective. 
Each ${\WH}_{(j)}$ is either a finite Coxeter group or an affine Weyl group \cite[V.3.9,VI.2.5]{bourbaki}. 
Accordingly, the fundamental chamber decomposes as $\AA=V_{\HH}+\prod_{j=1}^r\AA_{(j)}$ where each $\Abar_{(j)}$ is a fundamental domain for the action of ${\WH}_{(j)}$ on $E_{(j)}$ and it is either a simplex or a simplicial cone. Similarly, ${\mathcal P}(\HH)$ decomposes. Two faces $F=V_{\HH}+ \prod_{j=1}^rF_{(j)}$ and $F'=V_{\HH}+\prod_{j=1}^rF'_{(j)}$ in ${\mathcal P}(\HH)$ are separated by a single wall if and only if there is $j$ such that $F_{(j')}=F'_{(j')}$ for every $j'\neq j$ and $F_{(j)}$ and $F'_{(j)}$ are separated by a single wall in $\HH_{(j)}$. 

We will say that a flat $L$ {\em lies over $\AA$} if $L=|L\cap\Abar|$. This is the case if and only if $L$ is the intersection of some of the walls of $\AA$. 

\begin{remark}\label{rem:lie_over}A flat $L'$ in $E$ is always $\WH$-conjugate to a flat $L$ lying over $\Abar$. 
Indeed for any $\AA_{L'}\in{\mathcal C}(\HH^{L'})$ there is an ${\mathcal A}'\in{\mathcal C}(\HH)$ such that $\Abar_{L' }=\Abarp\cap L'$. Then, there is $w\in\WH$ such that $w\Abar_{L'}=\Abar\cap wL'$ and $L:=wL'$ lies over $\Abar$.
\end{remark}

We consider the group $\hWH:={\rm Stab}_{{\rm Aut}(E)}(\HH)$. It is the normaliser of $\WH$ in ${\rm Aut}(E)$ and  preserves ${\mathcal C}(\HH)$. 
%Let $A$ be the subgroup of ${\rm Aut}(E)$ stabilising $\overline{\AA}$. 
Then $\hWH\simeq {\rm Stab}_{\rm Aut(E)}(\AA)\ltimes \WH$, (proof as in \cite[V.2.3]{bourbaki}) and for any subgroup $\WW$ of $\hWH$ containing $\WH$ we have $\WH\lhd \WW$ and $K:=\WW/\WH\leq \hWH/\WH \simeq {\rm Stab}_{\rm Aut(E)}(\AA)$. Thus, all  subgroups of $\hWH$ containing $\WH$ are of the form $\WK=K\ltimes\WH$ for some $K\leq {\rm Stab}_{\rm Aut(E)}(\AA)$. Whenever we write $\WK$ we will always mean such a group.\label{WK}

\begin{remark}If $\WH$ is essential, then ${\rm Stab}_{\rm Aut(E)}(\AA)$ is finite. Indeed, it permutes the walls of $\Abar$ and therefore it permutes the (finitely-many) elements of the set ${\mathcal F}$ consisting of minimal dimensional faces of $\Abar$ that are not fixed by $\WH$. 
Let $\sigma\colon {\rm Stab}_{\rm Aut(E)}(\AA)\to {\rm Perm}({\mathcal F})$ be the corresponding group morphism. The elements in ${\mathcal F}$ are products of half lines and points and if $s\in {\rm ker}(\sigma)$, then $s$ must fix each of these faces pointwise because it is an Euclidean transformation. Hence $s=\id$ and ${\rm Stab}_{\rm Aut(E)}(\AA)$ is finite.
\end{remark}

 If the action of $\WH$ is not essential, then ${\rm Stab}_{\rm Aut(E)}(\AA)$ is never finite as it contains all translations by vectors in $V_{\HH}$. We will say that $K$ is admissible if its action  is obtained by pull-back of an action on the affine space $E^{ess}\simeq E/V_{\HH}$, i.e., if it satisfies $k(x+v)=kx+v$ for any $k\in K$, $x\in\Ec$ and $v\in V_{\HH}$. The $K$-action is trivial in the direction of $V_{\HH}$ and $K$ is always finite in this case. {\em From now on $K$ is assumed to be admissible.}  
%In this case, {\em we shall only consider extensions $\WK$ where the action of $K$ is obtained by pull-back of an action on the affine space $E^{ess}\simeq E/V_{\HH}$,} i.e., such that $k(x+v)=kx+v$ for any $k\in K$, $x\in\Ec$ and $v\in V_{\HH}$ and hence has trivial action in the direction of $V_{\HH}$. In particular, $K$ will always be finite. 
Observe that if $\WH$ is not irreducible, $K$ may permute the components of $E$.

We end this section by a simple observation that will be needed in the sequel.
\begin{remark}\label{rem:fundamental}If $kwp\in\Abar$ for some $p\in\Abar$, $k\in K$ and $w\in\WH$,  then $wp\in\Abar\cap\WH p=\{p\}$. Thus, $kwp=kp$.
\end{remark}

\subsubsection{Complexification}\label{subsec:complex}For any real affine space or vector space $A$, we will indicate by $A_{\mathbb C}$ its complexification. \label{complex} By abuse of terminology we will also call affine reflection hyperplanes the complexification of the affine hyperplanes in $\HH$ in $\Ec$. We set $\HHc:=\{H_{\mathbb C}~:~H\in\HH\}$. Note that $\Vc$ acts on $\Ec$ by translations and that any point in $\Ec$ can be written as $x+iy\in E+iV$. Such a point lies in a complexified affine hyperplane $H_{\mathbb C}$ if and only if $x\in H$ and $y\in V(H)$. The intersection of complex hyperplanes in $\HH_{\mathbb C}$ induces a stratification on $\Ec$ and $\hWH$ acts on $\Ec$ stabilising $\HH_{\mathbb C}$. Let $L$ be a flat in $E$. By the description of the complex affine hyperplanes, the (complex) flat $\Lc$ is the affine space containing a point in $L$ and having direction $V(L)_{\mathbb C}$. Its generic part will be denoted by $\mathring{\Lc}$. It follows that if $V_{\HH}=\{0\}$, then all $0$-dimensional flats lie in $E$. 

For any $D\subset \Ec$ and $w\in \hWH$ we shall set $D^w:=\{x\in D~:~wx=x\}$. When needed, we will denote a translation along a vector $v\in V_{\mathbb C}$ by $\tau_v$ and for any group $\WW$ acting on a set $D$ and any $p\in D$ we shall denote by $\WW^p$ \label{Wp} the stabilizer of $p$ in $\WW$. 

For $l\in\Ec$ we shall denote by $\HH_l$ \label{Hl} the subarrangement of $\HH$  consisting of hyperplanes whose complexification contains $l$. Then, $\WW_{\HH}^{l}=\WW_{\HH_l}$, i.e., the subgroup of $\WH$ generated by the reflections with respect to the hyperplanes in $\HH_l$: if $l\in E$ this is \cite[V.3.3]{bourbaki}, so $\HH_l$ is again admissible. If $l=x+iy\in\Ec$ with $x\in E$ and $y\in V$, then $w\in\WHl$ implies that $w\in\WW_\HH^x$ and $w$ acts linearly on $V$ and $iV$ fixing the direction $y$. In other words, $\WHl$ is generated by the reflections with respect to those affine hyperplanes containing $x$ and whose direction contains $y$, i.e., the hyperplanes in $\HH$ whose complexification contains $l$. The group $\WHl$ is finite by \cite[V.3.3 Proposition 2,V.3.6 Proposition 4]{bourbaki}.

A fundamental region $\Abarc$ \label{complex-fundamental} for the complexified action of $\WH$ on $\Ec$ is given by the set of points $x+iy$ such that $x\in \Abar$ and $y$ lies in the unique fundamental domain containing $\Abar$ for the action of the finite group $\WW_\HH^x$. 

\subsubsection{The main problems}\label{sec:basic-problems}
Let $\HH$ and $K$ be admissible and let $L$ be a flat. 
We set $X(\HH,K,L):=\WK \Lc/\WK$.\label{XL} If $K=1$ we will write $X(\HH,L)$ rather than $X(\HH,1,L)$. 
%When $\HH$ and $K$ are clear from the context we shall simply write $X(L)$. 
%By Remark \ref{rem:lie_over} we can always make sure that $L$ lies over $\Abar$.  

%If $K=1$ and we want to insist on this, we shall also denote it by $X_{\HH}:=\WH\Lc/\WH$. 
This quotient is well-defined at the level of analytic varieties (\cite[Satz 21, p. 186]{german}), i.e., $\WK$ acts on the analytic variety $\Ec$ and on the saturation $\WK \Lc$ and the quotient exists: its functions are the $\WK$-invariant analytic functions on $\WK \Lc$. By local finiteness of $\HH$ if $L'\subset L$ are strata for $\HH$, then $X(\HH,K,L')$ is closed in $X(\HH,K,L)$ as every point in $X(\HH,K,L)$ has a neighbourhood such that the intersection with $X(\HH,K,L')$ is closed. Hence the analytic variety $\Ec/\WK$ is stratified by the varieties $X(\HH,K,L)$, where $L$ ranges among all flats for $\HH$. We will parametrise strata in Section \ref{sec:sette}. 

In the present paper we shall address the following problems:
\begin{enumerate}
\item Provide $X(\HH,K,L)$ of the structure of an affine algebraic variety.
\item Determine when  $X(\HH,K,L)$ is unibranch, respectively normal, respectively smooth.
\end{enumerate} 
Question 1 is non-trivial only when $\WH$ is infinite. Also, when $K=1$ and $\WH$ is finite the normality and smoothness questions are answered in \cite{richardson,broer,DR}. Note that when $K=1$ the quotient $X(\HH,L)$ is the product of the quotients corresponding to the irreducible factors of $\WH$. Our main goal is to answer question 1 and the normality question  when $\WH$ is an affine Weyl group and $K$ is a natural abelian group related to $\WH$. 

\begin{remark}\label{rem:reductive}
If $V_{\HH}\neq\{0\}$ we can always reduce to the essential situation of $\HH^{ess}$ in $E^{ess}$. Indeed any flat $L\subset E$ is of the form $L'+V_{\HH}$ for a flat  $L'$ in $E^{ess}$ and $\Lc=L'_{\mathbb C}+(V_{\HH})_{\mathbb C}$. For our choice of $K$ we have
\begin{align*}
X(\HH,K,L)&=\WK \Lc/\WK=\WK (L'_{\mathbb C}+(V_{\HH})_{\mathbb C})/\WK\\
&\simeq\WK L'_{\mathbb C}/\WK\times (V_{\HH})_{\mathbb C}= X(\HH^{ess},K,L')\times(V_{\HH})_{\mathbb C}. 
\end{align*}
\end{remark}

\subsubsection{The affine Weyl group case}\label{rem:center} Assume $\WH=W_{aff}$ is an affine Weyl group acting effectively  on $E$. Then, for some point in $E$ which we can set as an origin $O$, there are a root system $\Phi$ with basis $\Delta=\{\alpha_1,\,\ldots,\,\alpha_\ell\}$, co-root lattice $Q^\vee={\mathbb Z}\Phi^\vee$, co-weight lattice $P^\vee$ and Weyl group $W=(W_{aff})^{O}$ such that $W_{aff}=W \ltimes Q^\vee$ and $E=Q^\vee\otimes_{\mathbb Z}{\mathbb R}$.
%For $\alpha\in \Phi^+$ and  $m\in{\mathbb Z}$ we set $H_{\alpha,m}:=\{h\in E_{\mathbb C}~|~\alpha(h)=m\}$ so all affine hyperplanes in $\HH$ are of this form, whereas the hyperplanes in $\HH_0$ are the $H_\alpha:=H_{\alpha,0}$. 
In this situation chambers are usually called alcoves. If $\Phi$ is irreducible then we choose as fundamental domain the closure of the fundamental alcove ${\mathcal A}$, which is the open simplex with vertices\label{vertices}  $x_0:=0$ and $x_1:=\omega_1^\vee/d_1,\,\ldots,\,x_\ell:=\omega_\ell^\vee/d_\ell$, where the $d_i$'s are the coefficients of the simple roots in the expression of the highest root $-\alpha_0$ and the $\omega_i^\vee$ are the fundamental co-weights.\label{coweights} For convenience we shall  set $\omega_0^\vee:=0$. 

A special family of groups of the form $\WK$ can be obtained by taking a lattice  $N^\vee$ satisfying $Q^\vee\subset N^\vee \subset P^\vee$. Translation by vectors in $N^\vee$ stabilises $\HH$, and $W$ stabilises $N^\vee$, so $\WH\leq \langle \WH, N^\vee\rangle=W\ltimes N^\vee\leq\hWH$, hence $W\ltimes N^\vee\simeq \WK$ for $K\simeq N^\vee/Q^\vee$. The group  $K$ acts on $\Abar$ as follows. Let $\nu\in N^\vee$ and let $x\in \Abar$. Then, $x+\nu\in E$ hence there is a unique $y\in \WH(x+\nu)\cap\Abar$. We set thus $\nu\cdot x=y$. Even though $x+\nu$ depends on the choice of the representative of the coset $\nu+Q^\vee$, the element $y$ does not. This procedure defines an action because $\WH$ normalises $N^\vee/Q^\vee$. 

%Observe that $Q^\vee\triangleleft\WK$ and in this situation $\WK/Q^\vee\simeq$

\subsection{Algebraic groups notation}Until otherwise stated  $G$ will denote a complex connected reductive algebraic group with Lie algebra $\gg$ and $T$ will be a fixed maximal torus in $G$ with Lie algebra $\hh$, whereas $W$ will be the Weyl group and $X_\alpha$ will denote the root subgroup corresponding to the root $\alpha$. If we insist that $G$ is semisimple and simply connected we shall write $G_{sc}$ and $T_{sc}$ instead of $G$ and $T$. The conjugation and adjoint action of $G$ on itself and $\gg$, respectively, will be denoted by a dot. The center of a group $C$ (a  Lie algebra $\cc$, respectively) will be denoted by $Z(C)$ ($\zz(\cc)$, respectively). 
The identity component of $H\leq G$ will be indicated by $H^\circ$.

If an algebraic group $H$ acts on a variety $Y$ we denote by $Y^{reg}$ the set of points in $Y$ whose $H$-orbit have maximum dimension. For $\psi$ an automorphism of a variety $Y$ we shall denote by $Y^\psi$ the set of points of $Y$ which are fixed by $\psi$. 

\subsection{Main motivation: Jordan classes and sheets in $G$}\label{sec:strata-Jordan}

The geometry of the stratifications induced by the decomposition of $\gg$ or $G$ into Jordan classes is the main motivation for our study.

A Jordan class in $\gg$ is an equivalence class with respect to the following equivalence relation: $x,y \in\gg$, with Jordan decomposition $x=x_s+x_n$ and $y=y_s+y_n$, respectively, are equivalent if, up to $G$-action, $x_s$ and $y_s$ have the same centraliser $\cc$ in $\gg$ and the nilpotent orbits represented by $x_n$ and $y_n$ in $\cc$ coincide, \cite{BK}. As a set, the class of $x$ is ${\JJ}(x)=G\cdot(\zz(\cc)^{reg}+x_n)$. 

A Jordan class in $G$ is  is an equivalence class with respect to the following equivalence relation: $x,\, y\in G$, with Jordan decomposition $x=su$ and $y=rv$, respectively,
are equivalent if, up to $G$-action, $s$ and $r$ have the same connected centraliser $M$ in $G$, $s\in Z(M)^\circ r$,  and the unipotent classes in $M$ represented by $u$ and $v$ coincide. As a set, the class of $x$ is ${J}(x)=G\cdot ((Z(M)^\circ s)^{reg}u)$. 

Jordan classes are finitely-many, locally closed, irreducible and smooth and the closure of a Jordan class is a union of Jordan classes \cite{lusztig-inventiones}, \cite[Proposition 4.9]{gio-espo}. Their closures form a stratification of $G$ and $\gg$, respectively. We consider the categorical quotient maps $G\to G/\!/G$ and $\gg\to \gg/\!/G$. The images of the strata of the respective Jordan stratifications form stratifications of $G/\!/G$ and $\gg/\!/G$ where the strata are of the form $\overline{J}/\!/G$ and $\overline{\JJ}/\!/G$, respectively, with $J$  and $\JJ$ semisimple Jordan classes, i.e., consisting of semisimple elements. We call them the Jordan stratifications of $\gg/\!/G$ and  $G/\!/G$. In the Lie algebra case, strata are closures of Luna strata \cite[III.2]{luna}, whereas in the group case, they are irreducible components of closures of Luna strata.

When $\WH$ is finite, we have the following correspondence, stemming from \cite[\S 2]{DR}.
\begin{proposition}\label{prop:Lie-stratification}Let $Q^\vee$ be the co-root lattice of $\gg$, let $E=Q^\vee\otimes_{\mathbb Z}{\mathbb R}$ and let $\HH$ be the  (finite) hyperplane arrangement of $W$. The stratification of $\Ec/\WH=\hh/W$ corresponds to Luna stratification of $\gg/\!/G$ through Chevalley isomorphism $\gg/\!/G\to \hh/W$.
\end{proposition}
\pf Let  $\JJ\subset\gg$ be a semisimple Jordan class. Then $\overline{\JJ}\cap \hh=W\zz$ where $\zz$ is the center of a standard Levi subalgebra $\cc$ containing $\hh$ and $\JJ$ is completely determined by $\overline{\JJ}\cap \hh$. Hence, $\zz$ is the intersection of the reflection hyperplanes corresponding to any choice of simple roots of $\cc$, \cite{DR}. Any such $\zz$ gives rise to a unique Jordan class and by construction, the closed sets $\overline{\JJ}/\!/G$ and $W\zz/W$ correspond through Chevalley isomorphism.  
\epf

We describe now a similar correspondence in the case of $\WH=W_{aff}=W\ltimes Q^\vee$ infinite, affine. Let $G_{sc}$ be the semisimple simply connected group with Weyl group $W$ and coroot lattice $Q^\vee$ and let $e=\exp(2\pi i -)\colon \hh\to T_{sc}$ be the exponential map. \label{exp} Its kernel is $Q^\vee$. This map realizes $\Ec:=\hh$ as the universal cover of $T_{sc}$ and $Q^\vee$ is the fundamental group of $T_{sc}$. Thus, the map $e$  induces an isomorphism of analytic varieties  
\begin{align}\label{eq:iso-sc}\Ec/\WH=\Ec/W_{aff}=\Ec/(Q^\vee\rtimes W)\simeq(\Ec/Q^\vee)/W\simeq T_{sc}/W.\end{align} 
More generally, let $G_{sc}\to G$ be a central isogeny inducing the natural projection $\pi\colon T_{sc}\to T$. Then $K:={\rm Ker}(\pi)\simeq N^\vee/Q^\vee$ for some lattice $Q^\vee\subset N^\vee \subset P^\vee$. Let $\HH$ be the arrangement of $W_{aff}=Q^\vee\rtimes W$. Then $W\ltimes N^\vee\simeq \WK= K\ltimes(W\ltimes Q^\vee)$ is one of the subgroups introduced in Subsection \ref{subsec:basic}. The $W$-action on $N^\vee$ induces the trivial action on $K=N^\vee/Q^\vee$ and $Q^\vee\triangleleft\WK$, so $\WK/Q^\vee\simeq W\times K$. Therefore, the map $e_\pi:=\pi\circ e$ induces an isomorphism of analytic varieties  
 \begin{align}\label{eq:iso-general}
 \Ec/\WK\simeq(\Ec/Q^\vee)/(W\times K)\simeq T_{sc}/(W\times K)\simeq (T_{sc}/K)/W\simeq T/W.
   \end{align}

\begin{proposition}\label{prop:stratifications}Let $G$, $\Phi$, $Q^\vee$, $K$ and $W_{aff}$ be as above and let $\HH$ be the arrangement of $\Ec=\hh$ such that $\WH=W_{aff}$. The isomorphism of analytic varieties $\hh/\WK\simeq G/\!/G$ given by \eqref{eq:iso-general} followed by Chevalley isomorphism identifies the stratification induced by $\HH$ and $K$ on $\Ec/\WK$ with the Jordan stratification in $G/\!/G$. 
\end{proposition}
\pf Chevalley isomorphism $T/W\simeq G/\!/G$ is induced by the inclusion $T\subset G$ and it was shown in the proof of \cite[Theorem 2]{proceedings-joseph} that if $J=G\cdot (Z(M)^\circ s)^{reg}$ for $s\in T$ and $M=G^{s\circ}$, then the stratum $\overline{J}/\!/G$ is identified with $W\cdot (Z(M)^\circ s)/W$. The possible pairs $(M, Z(M)^\circ s)$ run through the set of $W$-conjugacy classes of pairs of the form $(G^{r\circ},Z(G^{r\circ})^\circ r)$ for some $r\in T$.
%and $Z(M)^\circ r$ the coset in $Z(M)$ containing $r$. 
Let $L$ be a flat for $\HH$ and let $x\in \interior$, so $x$ does not lie in any hyperplane in $\HH$ not containing $L$. Thus, there is a subset $\Psi\subset\Phi^+$ such that $\Lc=x+\bigcap_{\alpha\in \Psi}{\rm Ker}\,\alpha$ and $\beta(x)\in{\mathbb Z}$ if and only if $\beta\in \Psi$. Then $M:=G^{e_\pi(x)\circ}=\langle T,\,X_{\pm\alpha}~:~\alpha\in \Psi\rangle$,  
$\bigcap_{\alpha\in \Psi}{\rm Ker}\,\alpha=\zz({\rm Lie}(M))$ hence 
$e_\pi(\Lc)=Z(M)^\circ e_\pi(x)$. Conversely, if $(M,Z(M)^\circ r)$ is a pair corresponding to a semisimple Jordan class, then $r=e_\pi(x)$ for some $x\in\hh$ and, taking the unique flat $L$ such that $x\in \interior$, the above argument shows that $e_\pi(\Lc)=Z(M)^\circ r$. The series of identifications in \eqref{eq:iso-general} maps then $X(\HH,K,L)=\WK\Lc/\WK$ to $W(Z(M)^\circ r)/W$.\epf
%Restricting $e_\pi$ to $e_\pi^{-1}(Z(M))=\{x\in \hh~:~\beta(x)\in\ZZ,\,\forall\beta\in\Pi\}$ gives again a covering with group $N^\vee$. The connected components  of $e_\pi^{-1}(Z(M))$ are complex flats for $\HH$. We  observe that connected components of a cover are  mapped onto connected components of the base and $Z(M)^\circ s$ is a connected component of $Z(M)$. Hence, $e_\pi^{-1}(Z(M)^\circ s)$ is a $N^\vee$-orbit of a flat $\Lc$ and $N^\vee \Lc/N^\vee\simeq Z(M)^\circ s$. Hence, $e_\pi$ identifies $\WK\Lc/\WK$ with $W (Z(M)^\circ s)/W$.

A similar statement for $G$ reductive could be proved, provided $K$ is allowed to be a discrete group containing suitable translations in the direction of vectors in $V_{\HH}$. 

The above interpretation of the stratification in $\Ec/\WH$ allows us to answer question 1 from Section \ref{sec:basic-problems}.

\begin{corollary} The stratum $X=X(\HH,K,L)$ is an affine algebraic variety for any admissible $\HH$ and $K$ and any flat $L$.
\end{corollary}
\pf  Since $X=K X(\HH,L)/K$ and $K$ is finite, it is enough to prove the statement for $X(\HH,L)$. By Remark \ref{rem:reductive} we may assume that $\HH$ is effective.  Without loss of generality we assume $\WH=W_{aff}=W\ltimes Q^\vee$, for $W$ the Weyl group of a complex semisimple algebraic group $G$ and $Q^\vee$ its coroot lattice.
%$X\simeq (K(W e(\Lc)/W))/K$, hence $X$ is an affine algebraic variety.
%Let $e(\Lc)=Z(M)^\circ s\subset T_{sc}$ as above. Then, 
By Proposition \ref{prop:stratifications} we have $X(\HH,L)\simeq \overline{J}/\!/G$, for some semisimple Jordan class $J$, hence it can be equipped with an affine algebraic variety structure.
%where $Z(M)^\circ s$ is a shifted torus in $T_{sc}$ and $W$ is a finite group acting on $T_{sc}$. 
\epf

Our approach to the problems we address will mainly rely on a local study and we will make use of the following key notion from \cite[1.7]{he}.

\begin{definition}Two pointed (algebraic or analytic) varieties $(Y_i,y_i)$, for $i=1,2$ are smoothly equivalent if there exist a pointed variety $(Y,y)$ and two smooth maps $\varphi_i\colon Y\to Y_i$ for $i=1,2$ such that $\varphi_i(y)=y_i$. 
\end{definition}
Smooth equivalence will be denoted by $\sim_{se}$. \label{se} Two pointed (algebraic or analytic) varieties $(Y_i,y_i)$ for $i=1,2$ satisfying $\dim Y_1=\dim Y_2+d$, are smoothly equivalent if and only if  $(Y_1,y_1)$ and $(Y_2\times{\mathbb A}^d,(y_2,0))$ are locally analytically isomorphic, \cite[Remark 2.1]{KP}. In particular, if $d=0$  there is a local analytic isomorphism between neighbourhoods of $y_1$ and $y_2$ mapping $y_1$ to $y_2$.  Smooth equivalence preserves the properties of being unibranch, normal, or smooth, \cite[Expos\'e XII, Proposition 2.1(vi), Proposition 3.1 (vii)]{SGA}.
 %\begin{proposition}\label{lem:analytic}
%Let $X$ be an algebraic variety and let $x\in X$. Then $X$ is unibranch, respectively normal, in $x$ if and only if the analytic variety $X^{an}$ associated with $X$ is unibranch, respectively normal, in $x$.
%\end{proposition}
%\pf This is \cite[Expos\'e XII, Proposition 2.1(vi), Proposition 3.1 (vii)]{SGA}.
%\epf
Thus, question 2 from Section \ref{sec:basic-hyperplanes}, for $\WH=W_{aff}$ and $K\leq P^\vee/Q^\vee$ translates into the following question:
\begin{center}{\em When is a stratum $\overline{J}/\!/G$  unibranch, respectively normal, respectively smooth?}
\end{center}

If $K=1$ and $\WH$ is a finite Weyl group, i.e., when the strata correspond to Jordan strata in $\gg/\!/G$,  it was shown in \cite{broer} that  $\overline{\JJ}/\!/G$ is normal if and only if it is smooth. We will show that for $K=1$ this is always the case.
%We will give an explanation of this phenomenon and show that it holds also when $\WH=W_a$, provided $K=1$. 

\begin{remark}
Let ${\mathcal G}\in\{G, \gg\}$, and let ${\mathcal S}$ be a sheet in ${\mathcal G}$, i.e., an irreducible component of the locally closed subset ${\mathcal G}_{(d)}$ consisting of the union of all the $G$-orbits in ${\mathcal G}$ of dimension $d$ for some $d\geq0$. By \cite{BK}, \cite[Propositions 5.1, 5.3]{gio-espo}, every sheet ${\mathcal S}$ contains a unique dense Jordan class. It was observed in \cite{broer}, \cite[\S 4]{proceedings-joseph} that the collection of quotients $\overline{\mathcal S}/\!/G$ where ${\mathcal S}$ runs among all sheets in ${\mathcal G} $ coincides with the collection of quotients of closures of semisimple Jordan classes. Hence, a complete list of normal or smooth strata in ${\mathcal G}/\!/G$ is also the complete list of normal or  smooth quotients of closures of sheets in ${\mathcal G}$. 
\end{remark}

\section{The normalisation of $X(\HH,K,L)$}

In this section we describe the normalisation of $X(\HH,K,L)$. 
%The symbol $L$ will always denote a flat for $\HH$. 
%for all choices of $\WH$, $K$ and $L$.
%Here, if $\WH$ is finite, it fixes a point that we take as the origin $O$, so $D(L)$ is an intersection of hyperplanes in $\HH_0$. 
By abuse of notation we will sometimes say that $X(\HH,K,L)$ is normal, respectively unibranch, respectively smooth, at $l\in \Ec$ if it is normal, respectively unibranch, respectively smooth, at the class $\overline{l}$ in $X(\HH,K,L)$ represented by $l\in L$. We set  \label{gamma}
\begin{align}\label{eq:gammaK}&\Gamma_{\HH,K,L}={\rm Stab}_{\WK}(L)={\rm Stab}_{\WK}(\Lc),\quad
\Gamma_{\HH,L}:={\rm Stab}_{\WH}(L)={\rm Stab}_{\WH}(\Lc).\end{align} 
Observe that $\Gamma_{\HH,L}$ preserves the components of $\Ec$ if $\WH$ is not irreducible. When $\WH$ is an affine Weyl group $\WH=W_{aff}=W\ltimes Q^\vee$ we will also need the group
\begin{align*}\Gamma_{W,L}:={\rm Stab}_W(e(\Lc)).\end{align*} 
We consider the quotients $\tilde{X}(\HH,L)=\Lc/\Gamma_{\HH,L}$ and $\tilde{X}(\HH,K,L):=\Lc/\Gamma_{\HH,K,L}$. If $\WH$ is finite then  $\tilde{X}(\HH,K,L)$ and $\tilde{X}(\HH,L)$ are normal affine varieties. We discuss the case of infinite components in $\WH$.

\begin{lemma}\label{lem:series}Assume $\WH=W_{aff}$ and let $L$ be a flat with $v\in L$. Then: 
\begin{enumerate}[label=(\roman*)]
 \item\label{item:isom} The natural projection of $\WH$ onto $W$ induces an isomorphism 
 \begin{equation*}\Gamma_{\HH,L}/(Q^\vee\cap V(L))\simeq \Gamma_{W,L}.\end{equation*}
 \item\label{item:isom2} The map $e$ induces an isomorphism from $\tilde{X}(\HH,L)$ to $e(\Lc)/\Gamma_{W,L}$.
 \item\label{item:normal} $\tilde{X}(\HH,L)$ is a normal algebraic variety.
\end{enumerate}
\end{lemma}
\pf  \ref{item:isom}. We consider the composition $\Gamma_{\HH,L}\subset \WH\to W$. The kernel 
%consists of those translations along vectors in $Q^\vee$ preserving $L$, i.e., vectors in  
is $Q^\vee\cap V(L)$. Assume $\tau\sigma \in \Gamma_{\HH,L}\leq Q^\vee\rtimes W$. Then, 
%for every $l\in\Lc$ we have 
$e(\Lc)=e(\tau\sigma(\Lc))=e(\sigma(\Lc))=\sigma(e(\Lc))$, so
%$\tau\sigma (l)\in\Lc$ so $\sigma(e(l))=e(\sigma(l))=e(\tau\sigma(l))\in e(\Lc)$, whence 
$\sigma\in \Gamma_{W,L}$.  We show surjectivity. If $\gamma\in \Gamma_{W,L},$  
then $\gamma(\Lc)$ is a connected component of $\Lc+Q^\vee$. Hence, $\gamma(\Lc)=v+q+V(L)$ for some $q\in Q^\vee$ and $\tau_{-q}\gamma\in  \Gamma_{\HH,L}$ is a pre-image of $\gamma$.
%, hence the image of the composition is $\Gamma_{W,L}$.
% \noindent\ref{item:isom2}. By \ref{item:isom}, we have $\tilde{X}(\HH,L)\simeq (\Lc/Q^\vee\cap V(L))/(\Gamma_{\HH,L}/Q^\vee\cap V(L))\simeq e(\Lc)/\Gamma_{W,L}$.

\noindent\ref{item:normal} follows from \ref{item:isom2} because $e(\Lc)$ is a smooth algebraic variety and $\Gamma_{W,L}$ is finite. 
\epf

For $L$ a flat and $l\in \Lc$ we set \label{omega}
\begin{equation*}
\Omega(L,\HH,K)_l:=\left\{w\Lc~:~ w\in\WK,\; l\in w\Lc\right\}.
\end{equation*}
Since $\HH$ is locally finite, the above set is finite. In addition, if $L'\subset L$ is the unique flat such that  $l\in\interiorp_{\mathbb C}$, then $\Omega(L,\HH,K)_l=\Omega(L,\HH,K)_{l'}$ for any $l'\in \interiorp_{\mathbb C}$. 

\begin{lemma}\label{lem:orbits}Let $l\in\Lc$. There is a natural bijection 
between $(\WK l\cap\Lc)/\Gamma_{\HH,K,L}$ and $\Omega(L,\HH,K)_l/\WKl$.
\end{lemma}
\pf  The map assigning to $wl\in\WK l\cap\Lc$ the $\WKl$-orbit $\WW_{\HH,K}^{l}w^{-1}\Lc$ in $\Omega(L,\HH,K)_l$ is surjective on $\Omega(L,\HH,K)_l/\WKl$ by construction. 
In addition, $w_1 l$ and $w_2l$ are mapped to the same $\WKl$-orbit if and only if $\omega w_1^{-1}\Lc=w_2^{-1}\Lc$ for some $\omega\in \WKl$, i.e., if and only if  $w_2l=w_2\omega l \in\Gamma_{\HH,K,L}w_1l$.
% This happens if and only if $w_2l\in \Gamma_K w_1 l$.  
The defined map induces the desired bijection.
\epf

\begin{proposition}\label{prop:normalization}The quotient $\tilde{X}(\HH,K,L)$ is the normalisation of $X(\HH,K,L)$. 
\end{proposition}
\pf Assume first $K=1$. By Remark \ref{rem:reductive} we can reduce to the case that $\WH$ acts effectively and $\WH$ is either a finite  Coxeter group or an affine Weyl group.  If $\WH$ is a finite Weyl group the statement is \cite[Korollar 6.4 (b)]{bo}. If $\WH$ is finite, then 
the composition of the closed embedding $\Lc\to \WH\Lc$ with the quotient map $\WH\Lc\to\WH\Lc/\WH$ is a surjective finite morphism and it factors through $\tilde{X}(\HH,L)$.  We prove that the (finite) induced morphism $\tilde{X}(\HH,L)\to X(\HH,L)$ is generically injective. If  $l=x+iy\in\Lc$  and $x$ lies in a chamber $C$ of $\HH^L$, then the only hyperplanes of $\HH$ containing $l$ are those containing $\Lc$. If for such an $l$ we have $wl\in\Lc$ for some $w\in\WH$, then $w^{-1} \Lc\subset \bigcap_{\genfrac{}{}{0pt}{}{H\in\HH}{l\in H}}H=\Lc$, hence $w^{-1}\Lc=\Lc$, so $\WH l\cap \Lc=\Gamma_{\HH,L} l$. Thus, the morphism is generically bijective and it is the normalisation map.  If $\WH$ is an affine Weyl group, the result is obtained combining Proposition \ref{prop:stratifications}, and Lemma \ref{lem:series} (b) with \cite[Theorem 2]{proceedings-joseph}. 

Let now $K$ be arbitrary. The composition $\Gamma_{\HH,K,L}\to\WK\to \WK/\WH=K$ has kernel $\Gamma_{\HH,L}$, hence $\Gamma_{\HH,K,L}/\Gamma_{\HH,L}$ is finite and $\tilde{X}(\HH,K,L)$ is the quotient of the normal algebraic variety $\tilde{X}(\HH,L)$ by the action of this finite group, thus it is normal. 
The composition of the normalisation map $\tilde{X}(\HH,L)\to X(\HH,L)$ with the closed embedding $X(\HH,L)\to KX(\HH,L)$ and the quotient map $KX(\HH,L)\to KX(\HH,L)/K=X(\HH,K,L)$ is a finite surjective morphism factoring through  $\tilde{X}(\HH,K,L)$. The previous argument shows that the induced map $\tilde{X}(\HH,K,L)\to {X}(\HH,K,L)$ is generically injective whence it is  the normalisation map.
 \epf

\begin{corollary}\label{cor:condition_unibranch} The following three conditions on $l\in\Lc$ and $X(\HH,K,L)$ are equivalent:
\begin{enumerate}[label=(\roman*)]
\item\label{item:uni}The variety  $X(\HH,K,L)$ is unibranch at $l$;
\item\label{item:orbit_gamma} $\WK l\cap\Lc=\Gamma_{\HH,K,L} l$;
\item\label{item:orbit}
%\begin{equation*}
$\Omega(\HH,K,L)_l=\{w\Lc~:~w\in \WKl\}$.
% %.\end{equation*}o
\end{enumerate}
\end{corollary}
\pf The equivalence of \ref{item:uni} and \ref{item:orbit_gamma} follows from Proposition \ref{prop:normalization}. The equivalence between  \ref{item:orbit_gamma}  and \ref{item:orbit} is a consequence of Lemma \ref{lem:orbits}.
\epf

The following proposition is the analogue of \cite[Theorem A]{richardson}, which corresponds to the case $K=1$, $\WH=W$, a finite Weyl group.  When $\WH$ is an affine Weyl group, a similar statement follows 
combining Proposition \ref{prop:stratifications}, Lemma \ref{lem:series} (ii) and \cite[Theorem 2 (b)]{proceedings-joseph}.
%
%\pf The variety $X(\HH,K,L)$ is normal if and only if $X(\HH,K,L)\simeq\tilde{X}(\HH,K,L)$. Thus, if $X(\HH,K,L)$ is normal, the composition of maps $\tilde{X}(\HH,K,L)\to X(\HH,K,L)\to \Ec/\WK$ is a closed immersion, i.e., $\tilde{X}(\HH,K,L)=\Lc/\Gamma_{\HH,K,L}$ is a closed subvariety of $\Ec/\WK$ and therefore the map is surjective. Conversely, the map factors through ${\mathbb C}[\WK\Lc]^{\WK}={\mathbb C}[X(\HH,K,L)]$ \comgio{??} and if it is surjective, then $\tilde{X}(\HH,K,L)=\Lc/\Gamma_{\HH,K,L}$ is a closed subscheme of $X(\HH,K,L)$. \comgio{????}However, the only ideal of ${\mathbb C}[X(\HH,K,L)]$ whose zero locus is $X(\HH,K,L)$ is the $0$ ideal so $\tilde{X}(\HH,K,L)\simeq X(\HH,K,L)$\comgio{???}. 
%\epf
%The following proposition is the analogue of \cite[Theorem A]{richardson}, which corresponds to the case $K=1$, $\WH=W$, a finite Weyl group. The case of $\WH$ an affine Weyl group, and $K\simeq P^\vee/Q^\vee$ is covered by \cite[Theorem 2]{proceedings-joseph}. 

\begin{proposition}\label{prop:richardson}Assume $\WH$ is finite. The variety $X=X(\HH,K,L)$ is normal if and only if the map 
\begin{align*}{\mathbb C}[\Ec]^{\WK}\to {\mathbb C}[\Lc]^{\Gamma_{\HH,K,L}}\end{align*} induced from the natural restriction map is surjective.
\end{proposition}
\pf We consider the composition of the normalisation map $\tilde{X}(\HH,K,L)\to X$ with the closed immersion $X\to \Ec/\WK$ and the corresponding maps of algebras of regular functions. Then $X$ is normal if and only if $\tilde{X}(\HH,K,L)\to X$ is an isomorphism, which is equivalent to surjectivity of the algebra map ${\mathbb C}[\Ec]^{\WK}\to {\mathbb C}[\WK\Lc]^{\WK}$ whence 
to surjectivity of the restriction map ${\mathbb C}[\Ec]^{\WK}\to {\mathbb C}[\Lc]^{\Gamma_{\HH,K,L}}$ by diagram chasing. \epf

\section{Local geometry of strata}

In this section we begin our local study of strata $X(\HH,K,L)$ around the class $\overline{l}$, for $l\in L$. %Observe that we can always choose $L$ and a representative $l$ of $\overline{l}$ in $\Lc$ such that $l\in\Abarc$. 

We will show that normality and smoothness of a stratum can be checked in special points in the minimal strata contained in $X(\HH,K,L)$. In order to do so, we will study the hyperplane arrangements $\HH_l$ for $l\in L$. Since $\HH_l$ is admissible, $\WHl$ permutes simply transitively the chambers in  
${\mathcal C}(\HH_l)$. In addition, if $\AA'\in {\mathcal C}(\HH)$ and $l\in\overline{\AA'}$, then $\AA'\subset C$ for a unique $C\in{\mathcal C}(\HH_l)$. 
%Conversely, the walls of any $C\in {\mathcal C}(\HH_l)$ are a subset of the walls of a unique chamber $\AA'$ containing $l$ in its closure. This sets a bijective correspondence between the set of chambers in $\HH_l$ and the set of chambers of $\HH$ containing $l$ in their closures. 

\begin{lemma}\label{lem:finite}Let $l=x+iy\in \Abarc$. Then, $\WKl=K^l\ltimes \WHl$. Hence $\WW_{\HH,K}^z$ is finite for any $z\in \Ec$.  
\end{lemma}
\pf We need to prove $\subseteq$. Assume first that $l\in \Abar \subset E$. Let $kw\in \WKl$ with $k\in K$ and $w\in \WH$. Remark \ref{rem:fundamental} gives $l=kwl=kl$ whence $k\in K^l$ and $w\in \WHl$. Thus, $\WKl=K^l\ltimes \WHl$. 

Assume now $l=x+iy$ for $x\in \Abar$ and $y\in V$. Since $\WK$ preserves $E$, we have $\WKl\subset \WW_{\HH,K}^x\simeq K^x\ltimes\WW_{\HH}^x$, with $\WW_{\HH}^x=\WW_{\HH_x}$. Observe that, fixing $x$ as an origin of $E$ and identifying $E$ with $V$, the action of $K^x\ltimes\WW_\HH^x$ on $V$ is linear, so $\WKl=(K^x\ltimes\WW_\HH^x)^y$.

If $x\in {\mathcal A}$, then $\WW_\HH^x=\WHl=1$ and $\WKl=K^l$ and we are done. If, instead, $x\in\Abar\setminus{\mathcal A}$, then $y$ lies in the unique fundamental domain $D$ for the action of $\WW_\HH^x$ with $\Abar\subset D$. Since $K^x$ preserves ${\mathcal C}(\HH_x)$ and $\Abar$, it preserves $D$. Remark \ref{rem:fundamental} applied to $y$, $D$ and the group $K^x\ltimes\WW_\HH^x$ gives $(K^x\ltimes\WW_\HH^x)^y=K^l\ltimes \WHl$. Finally, $\WW_{\HH,K}^z$ is $\WH$-conjugate to $\WKl$ for some $l\in\Abarc$, hence it is finite 
because $K^l$ and $\WHl$ are so. 
\epf

For any flat $L$ we consider the groups:\label{WHKL}
\begin{align*}\WW_\HH^{L}:=\cap_{l\in\Lc}\WW_\HH^{l},&&\WW_{\HH,K}^{L}:=\cap_{l\in\Lc}\WW_{\HH,K}^{l},&&K^{L}:=\cap_{l\in\Lc}{K}^{l}\end{align*}
and the subset
\begin{align}\label{eq:ul}U_L:=\{l\in\Lc~:~\WKl=\WW_{\HH,K}^{L}\}.
 \end{align}

\begin{lemma}\label{lem:UL}
Let $L$ be a flat in $E$. The subset $U_L$ is a non-empty open subset of $\Lc$ contained in $\mathring{\Lc}$ and having non-empty intersection with $E$. By construction $U_L$ is the set of points in $L$ with minimum stabiliser in $\WK$. In addition, if $L$ lies over $\Abar$, then $\WW_{\HH,K}^L=K^{L}\ltimes \WW_{\HH}^{L}$ and $U_L\cap\Abar\neq\emptyset$.
\end{lemma}
\pf Assume first that $L$ lies over $\Abar$. Observe that $\Lc\cap\Abarc$ generates $\Lc$ as an affine space, hence 
\begin{align*}\WW_\HH^{L}=\cap_{l\in\Lc\cap\Abarc}\WHl,&&\WW_{\HH,K}^{L}=\cap_{l\in\Lc\cap\Abarc}\WKl,&&K^{L}=\cap_{l\in\Lc\cap\Abarc}{K}^{l}.\end{align*}
Clearly, $ K^{L}\ltimes \WW_\HH^{L}\leq \WW_{\HH,K}^{L}$. On the other hand, 
 if $kw\in\WW_{\HH,K}^{L}$ with $k\in K$ and $w\in\WH$, then  $k\in K^l$ and $w\in \WHl$ for any $l\in\Abarc\cap\Lc$ by Lemma \ref{lem:finite}, i.e., $\WW_{\HH,K}^{L}= K^{L}\ltimes \WW_\HH^{L}$.
By construction we have $\WW_{\HH,K}^{L}\leq \WKl$ for any $l\in\Lc$  and equality holds if and only if  $\WW_\HH^{L}=\WHl$ and $K^{L}=K^l$. 
The first condition holds if and only if $l$ lies in $\mathring{\Lc}$. The second one  holds if and only if $l\in U_1=\Lc\setminus \bigcup_{k\in K\setminus K^{L}}\Lc^k$. Thus, 
$U_L$ is the non-empty open set $\interiorc\cap U_1$. Also, $\interior$ is open in $L$ and $D=\bigcup_{k\in K\setminus K^{L}}L^k$ is a proper closed subset in $L$ and $\emptyset\neq\interior\cap (L\setminus D)\subset U_L$. So, $U_L\cap E\neq\emptyset$ and $\Abar\cap D\neq\emptyset$ by dimensional reasons.

Assume now that $L$ does not lie over $\Abar$. By Remark \ref{rem:lie_over} there always is $w\in \WH$ such that $wL$ lies over $\Abar$ and we have $U_L=w^{-1}(U_{wL})$.
\epf

For $l\in \Lc$,  we set:\label{Y}
%let $L'=\cap_{H\in\HH_l} H\subset \Lc$. We set:
\begin{equation}\label{eq:Y}Y(\HH,K,L)=\WK \Lc,\quad Y(\HH,K,L)_l=\bigcup_{\genfrac{}{}{0pt}{}{w\in \WK} 
{l\in w\Lc}}w\Lc.\end{equation}
Since $K$ preserves $\HH$, the set $Y(\HH,K,L)_l$ depends only on the flat $L'$ such that $l\in\interiorp$. The finite group $\WKl$ acts on $Y(\HH,K,L)$ and on $Y(\HH,K,L)_l$. 

%\begin{lemma}\label{lem:local}
%Let $l\in \Lc$. With the above notation, the analytic varieties $Y^{an}$ and $Y_l^{an}$ are locally isomorphic in  a neighbourhood of $l$. 
%The isomorphic neighbourhoods of $l$ can be chosen to be $\WKl$-stable and the isomorphism to be $\WKl$-equivariant.
%\end{lemma}
%\pf 

We recall a basic result.
\begin{lemma}(\cite[Anhang zu K. 7, Satz 21]{german}\label{lem:invariants}
Let $\WW$ be a discrete group acting properly discontinuously on an analytic variety $Y^{an}$ and let $X^{an}$ be the quotient of $Y^{an}$ by  $\WW$. If $\tilde{x}\in Y^{an}$ is mapped to $x$ in $X^{an}$ through the canonical quotient map and $H$ is the stabilizer of $\tilde{x}$ in ${\mathcal W}$, then there exists a small enough $H$-stable neighbourhood $U$ of $\tilde x$ in $Y^{an}$ such that $U/H$ can be identified with a neighbourhood of $x$ in $X^{an}$. \hfill$\Box$
\end{lemma}

\begin{proposition}\label{prop:local}Let $L$ be a flat, let $l\in\Lc$ and let 
$\overline{l}$ be the class of $l$. Then,
%$\overline{l}\in X(\HH,K,L)$, let $l\in \Lc\cap\Abar$ be a representative of $l$. Let $\overline{l}$ denote also the class of $l$ in $Y(\HH,K,L)_l/\WKl$. Then 
\begin{align}
\left(X(\HH,K,L),\overline{l}\right)\sim_{se} \left(Y(\HH,K,L)_l/\WKl,\,\overline{l}\right)
\end{align}
\end{proposition}
\pf By Lemma \ref{lem:invariants} with $X^{an}=X(\HH,K,L)$, $\tilde{x}=l$, $Y=Y(\HH,K,L)$ we have
\begin{align}
\left(X(\HH,K,L),\overline{l}\right)\sim_{se} \left(U/\WKl,\,\overline{l}\right)
\end{align}
where $U$ is a $\WKl$-invariant open neighbourhood of $l$ in $Y(\HH,K,L)$. Let $A$ be an analytic open neighbourhood of $l$ in $\Ec$ whose closure is compact. By local finiteness of $\HH$, the intersection $Y(\HH,K,L)\cap A$ is a finite union of translates of $\Lc$ and $A\cap Y(\HH,K,L)_l$ is the union of all irreducible components containing $l$ in the analytic variety $A\cap Y(\HH,K,L)$. Possibly reducing $A$ and invoking the finiteness result in Lemma \ref{lem:finite} we can make sure that 
$A$ is $\WKl$-stable and that $Y(\HH,K,L)\cap A=Y(\HH,K,L)_l\cap A$. 
Therefore  
\begin{align*}\left(U/\WKl,\,\overline{l}\right)\sim_{se}\left((U\cap A)/\WKl,\,\overline{l}\right)\sim_{se}\left(Y(\HH,K,L)_l/\WKl,\,\overline{l}\right)\end{align*}
concluding the proof. \epf
%
%Since $\HH$ is locally finite, the components $w\Lc$ occurring in  $Y(\HH,K,L)_l$ are finitely many and they are the irreducible components of  $Y(\HH,K,L)_l$. More precisely, $Y(\HH,K,L)_l$ is the union of all irreducible components of the analytic variety $Y(\HH,K,L)$ containing $l$
%so the two varieties are analytically isomorphic on a neighbourhood $U$ of $l$ in  $Y^{an}$. 
%$(Y(\HH,K,L),l)\sim_{se}(Y(\HH,K,L)_l)$.  
%By Lemma \ref{lem:finite} the group $\WKl$ is finite so we can always replace $U$ by the $\WKl$-stable open neighbourhood $\cap_{w\in\WKl}w(U)$ to ensure that the isomorphism is also $\WKl$-equivariant.
%\epf

%\begin{proposition}\label{prop:local} For any $l\in\Lc$, the analytic variety $X$ is locally isomorphic to $Y_l/\WKl$ around the point corresponding to $l$.
%\end{proposition}
%\pf  By Lemma \ref{lem:local}, $Y=\WK \Lc$ in a neighbourhood  $U$ of $l$ is isomorphic to $Y_l$, and the neighbourhood can be chosen to be $\WKl$-invariant. By Lemma \ref{lem:invariants}, $U/\WKl$ is isomorphic to a neighbourhood of both $X$ and $Y_l$ around the point corresponding to $l$.
%\epf

If $l\in\Abarc\cap \Lc$, then $\WK\simeq K^l\rtimes\WW_{\HH_l}$, with $\HH_l$ admissible. If $C$ is the unique chamber for $\HH_l$ containing $\AA$, then $C$ is a fundamental domain  for the action of $\WW_{\HH_l}$ and $K^l$ preserves $C$. In general if $\HH_l$ is not effective, $K^l$ is not necessarily admissible for $\HH_l$. However, if $l\in U_{L'}$ for some flat $L'$, then for any $H\in \HH$ we have $H\supset L'$ if and only if $H\in \HH_l$, so $V_{\HH_l}=\bigcap_{H\supset L'}V(H)=V(L')$ and $K^l=K^{L'}$ is admissible because it cannot contain translations in the direction of $V_{\HH_l}$.  In this case $\WKl\Lc/\WKl=X(\HH_l,K^l,L)$ with $\HH_l$ a {\em finite} arrangement.  We call it the {\em finite  counterpart} of $X(\HH,K,L)$ at $l$.

\begin{corollary}\label{cor:X'}%\comgio{vedere se serve anche per $l$ piu` generico, ed in caso aggiungere un'osservazione che l'affermazione e' vera anche per $l$ meno generico.}
Let $L'\subset L$ be flats,  let $l\in U_{L'}\cap\Abarc$ and $\overline{l}$ be its class in $X=X(\HH,K,L)$. Then $X$ is unibranch at $\overline{l}$ if and only if 
$\left(X,\overline{l}\right)\sim_{se}\left(X(\HH_l,K^l,L),\overline{l}\right)$. 
\end{corollary}
\pf If $X$ is unibranch at $\overline{l}$ then Proposition \ref{prop:local} and Corollary \ref{cor:condition_unibranch} give the desired equivalence. Conversely, if $\left(X,\overline{l}\right)\sim_{se}\left(\WKl \Lc/\WKl,\overline{l}\right)$, then $X$ is unibranch at $\overline{l}$ because the normalisation map $\tilde{X}(\HH_l,K^l,L)\to X(\HH_l,K^l,L)$ is bijective at $\overline{l}$.
\epf

Next Proposition will show that in order to check unibranchedness or normality of $X(\HH,K,L)$ at points in a flat $\Lc'\subset \Lc$, it will be enough to check it at one point in $U_{L'}$. It will allows us to reduce the verification of normality to suitable points in $\Abar$.  

\begin{proposition}\label{prop:constant}Let $L'\subset L$ be flats for $\HH$ % let $U_{L'}$ be as in \eqref{eq:ul}, 
and let $X=X(\HH,L,K)$.
%$X(\HH,K,L)\$ and let $\Lc'\subset\Lc$ be a flat in $\Ec$ and let $U_{L'}$ be as in \eqref{eq:ul}. 
Then:
\begin{enumerate}[label=(\roman*)]
\item\label{item:constant} $\left(X,\,\overline{l}\right)\sim_{se} \left(X,\,\overline{l'}\right)$ for any $l,l'\in U_{L'}$.
\item\label{item:uni_nor} If $X$ is unibranch, respectively normal, at some $l\in U_{L'}$, then it is again so at all $l'\in\interiorcp$. 
\end{enumerate}
\end{proposition}
\pf Statement \ref{item:constant} follows from Proposition \ref{prop:local} because $Y(\HH,K,L)_l=Y(\HH,K,L)_{l'}$ and $\WKl=\WW_{\HH,K}^{l'}$ and $l$ and $l'$ differ by a translation in $V(L)_{\mathbb C}$ which commutes with the $\WKl$-action. We prove  \ref{item:uni_nor}. Up to replacing $l$ and $\Lc$ by a point and a flat in their respective $\WH$-orbits, we may assume that $l\in \Abarc\cap U_{L'}$ and that $L'$ lies over $\Abar$. Assume first that $X$ is unibranch at $l$. Then, for every $l'\in\interiorcp$  we have 
\begin{align*}
\Omega(\HH,K,L)_{l'}&=\Omega(\HH,K,L)_l=\{w\Lc~:~w\in \WKl\}\\
&=\{w\Lc~:~w\in \WW_{\HH,K}^{L'}\}
%\left\{w\Lc~:~ w\in\WK,\; l'\in w\Lc\right\}&=\left\{w\Lc~:~ w\in\WK,\; l\in w\Lc\right\}\\
%&=\WKl\Lc={\WK}_{L'}\Lc\\
\subseteq \{w\Lc~:~w\in\WW_{\HH,K}^{l'}\}\subseteq \Omega(\HH,K,L)_{l'}
%\left\{w\Lc~:~ w\in\WK,\; l'\in w\Lc\right\}.
\end{align*}
where we applied Corollary \ref{cor:condition_unibranch} (iii) and Lemma \ref{lem:UL}. Thus, we have equality everywhere and $X$ is unibranch at $l'$.

Assume now that $X$ is normal at $l\in\Abarc\cap U_{L'}$ and let $l'\in\interiorcp$. We claim that $\WW_{\HH,K}^{L'}\lhd \WW_{\HH,K}^{l'}$. Indeed, if $kw\in\WW_{\HH,K}^{l'}$, then $kw\Lc'=\Lc'$ because $\Lc'$ is the only flat containing $l'$ in its generic part. Therefore for any $k'w'\in \WW_{\HH,K}^{L'}$ and any $l''\in \Lc'$ we have $kwk'w'((kw)^{-1}l'')=kw(kw)^{-1}l''=l''$. The quotient $H=\WW_{\HH,K}^{l'}/\WW_{\HH,K}^{L'}$ is finite by Lemma \ref{lem:finite}.  By the above discussion $X$ is unibranch at $l'$ and  Proposition \ref{prop:local} gives 
$\left(X,\overline{l'}\right)\sim_{se}\left(Y(\HH,K,L)_{l'}/\WW_{\HH,K}^{l'},\,\overline{l'}\right)=\left(Y(\HH,K,L)_{l}/\WW_{\HH,K}^{l'},\,\overline{l'}\right)=\left(X(\HH_l,K^l,L)/H,\overline{l'}\right)$.
% where we have used that $Y(\HH,K,L)_l=Y(\HH,K,L)_{l'}$
Hence it is enough to show that $X(\HH_l,K^l,L)$ is normal at $\overline{l'}$. This is the case because Corollary \ref{cor:X'} gives $\left(X,\overline{l}\right)\sim_{se}\left(X(\HH_l,K^l,L),\overline{l}\right)\sim_{se}
\left(X(\HH_l,K^l,L),\overline{l'}\right)$, where the second equivalence is induced by the translation along $l'-l\in V(L)_{\mathbb C}$ which commutes with the $\WKl$-action and maps a neighbourhood of $l$ in  $\interiorp$ to a neighbourhood  of $l'$ therein. 
\epf

\begin{corollary}\label{cor:dim0}The variety $X=X(\HH,K,L)$ is normal if and only if it is normal at all points in all strata $X(\HH,K,L')$ with $L'\subset L$ of minimal dimension. \end{corollary}
\pf If $X$ is not normal at some point in $\interiorcp\subset \Lc$, by Proposition \ref{prop:constant} \ref{item:uni_nor} it is not normal at any point in $U_{L'}$. As the non-normality locus is closed, $X$ is not normal at any point in $\overline{\WK U_{L'}/\WK}=X(\HH,K,L')$. If, instead $X$ is normal at all points in $\interiorcp\subset \Lc$, then the non-normality locus  may intersect $X(\HH,K,L')$ only at points in strictly lower dimensional strata  $X(\HH,K,L'')$ for $L''\subsetneq L'$. Hence, if non-empty, the non-normality locus must contain some minimal stratum contained in $X$. \epf

\begin{remark}\label{rem:real}\begin{enumerate}[label=(\roman*)]
\item The proof of Corollary \ref{cor:dim0} applied to $L'=L$ shows that $X(\HH,K,L)$ is always normal, hence unibranch, at every point in $\interiorc$.
\item If $L'\subset L$ is minimal, then $\interiorcp=\Lc'$ because there are no flats contained in $L'$ and $V(L')=V_{\HH}$. Also, if $k\in K$ fixes $l\in\Lc'$, then it fixes  $\Lc'=l+V(L')_{\mathbb C}$ so $K^{L'}=K^l$ and $U_{L'}=L'$.
By Proposition \ref{prop:constant}, normality or unibranchedness at points in a minimal stratum $X(\HH,K,L')\subset X(\HH,K,L)$ are thus guaranteed by normality or unibranchedness at  a real point in $L'$ or even in $L'\cap\Abar$ if $L'$ is chosen in its $\WK$-orbit so that this intersection is non-empty.  

\item If $\HH$ is effective, then minimal strata contained in $X(\HH,K,L)$ correspond to $0$-dimensional flats $L'$, which are all contained in $E$. Each point therein has a representative in $\Abar$. By Remark \ref{rem:reductive} one can always reduce to this case. 

\item If $\HH$ is finite, then it contains a unique minimal stratum $L'$ and for $l\in U_{L'}$ the finite counterpart $X(\HH_l,K^l,L)$ is normal if and only if it is normal at $\overline{l}$.  
\end{enumerate}
\end{remark}

 %\begin{equation}Y(\tilde{x})/W_{a,\tilde x}=W_{a,\tilde x}(v+\zz)/W_{a,\tilde x}=(W_{e(\tilde x),a})_{\tilde x}(v+\zz)/(W_{s,a})_{\tilde x}
%\end{equation}
%where we used Lemma \ref{lem:Y(x)} and unibranchedness for the first and the third equality and Remark \ref{rem:affinization} for the second one. 
%Observe that the proofs of Lemmata \ref{lem:Y(x)} and \ref{lem:local} do not make use of the simply connected condition on $G$ and can be applied to a more general setup.
%So, if we set $Y'({\tilde{x}})=\bigcup_{w\in W_{e(\tilde x),a\atop\tilde{x}\in w(v+\zz)} } w(v+\zz)$, since $W_{a,\tilde x}(v+\zz)/W_{a,\tilde x}$ is unibranch at the orbit of $\tilde x$
%we have  $W_{a,\tilde x}(v+\zz)/W_{a,\tilde x}=Y'(\tilde x)/(W_{e(\tilde x),a})_{\tilde x}$. By the proof of Lemma \ref{lem:local}, we have the statement.\epf

%\begin{remark}\label{rem:reductive}Assume $\WH$ is not essential on $E$ and $\Ec$. This happens if and only if there exists $\hh_0\subset H$ for every $H\in \HH_0$ \comgio{verfiicare}. Then, $\hh_0$ is contained in all strata and $\WH\Lc\simeq \hh_0+\WH(\Lc/\hh_0)$, so $\WH\Lc/\WH\simeq \hh_0+ \WH(\Lc/\hh_0)/\WH$  and the geometry of such a stratum can be reduced to the geometry of the stratum in $(\Ec/\hh_0)/\WH$ corresponding to $\Lc/\hh_0$. \end{remark}

Combining the results obtained so far we get the following characterisation of normality of a stratum.
\begin{theorem}\label{thm:isolated1}
A stratum $X=X(\HH,K,L)$ in $\Ec/\WK$ is normal if and only if the following two conditions hold:
\begin{enumerate}[label=(\roman*)]
\item\label{item:minimal_uni} $X$ is unibranch at a point in every minimal stratum it contains.
\item\label{item:minimal_normal}For any minimal stratum $X(\HH,K,L')\subset X$ there is $l\in \WK L'\cap\Abar$ such that $X(\HH_l,K^l,L)$ is normal. 
%$X(\HH_l,K^l,L)$ at points in minimal strata are all normal.
\end{enumerate}
\end{theorem}
\pf By Corollary \ref{cor:dim0} and Remark \ref{rem:real} the stratum $X$ is normal if and only if it is normal at a point $\overline{l}$ in  each minimal stratum, with representative $l\in\Abar$.

If $X$ is normal, then condition \ref{item:minimal_uni} holds. By Remark \ref{rem:real}, Corollary \ref{cor:X'} applies at all such $l$, hence we have normality of  $X(\HH_l,K^l,L)$ for any (or for an) $\overline{l}$ in each minimal stratum. 

Conversely, if  condition \ref{item:minimal_uni} holds, then  $X$ is unibranch at all points in every minimal stratum it contains by Remark \ref{rem:real} and  $\left(X,\,\overline{l}\right)\sim_{se} \left(X(\HH_l,K^l,L),\,\overline{l}\right)$ for all $\overline{l}$ in a minimal stratum by Corollary \ref{cor:X'}. If in addition condition \ref{item:minimal_normal} holds, then $X$ is normal at all such $\overline{l}$. We conclude by using Proposition \ref{prop:constant} \ref{item:constant}. 
\epf

\section{Necessary conditions for normality of strata}

In this Section we will provide necessary conditions to verify normality of a stratum in special cases. 

\subsection{The case $\dim \Lc=1$ and $\HH$ finite}
In this subsection $\HH$ is finite. Since $\HH$ is admissible, the chambers are finitely many, hence $\WH$ is a finite Coxeter group and by \cite[V.3.6 Proposition 4]{bourbaki} it fixes a point $p\in \cap_{H\in\HH}H\subset E$ and contains no translation.  Taking $p$ as an origin we identify $\Ec$ with $V_{\mathbb C}$ and similarly a flat $\Lc$ with its direction $V(\Lc)$. 

Assume $\dim\Lc=1$. Then, either $L=p+V_{\HH}=E^{\WH}$, or else $V_{\HH}=\{0\}$  and $\Ec^{\WH}=\{p\}\subset \Lc$. In the first case $\Lc$ is contained in all fundamental domains for the $\WH$-action and it is fixed pointwise by any admissible $K$, so $\WK=\Gamma_{\HH,K,L}$ and $L=E^{\WK}$. In the latter case, $p$ is the only point contained in any fundamental domain for $\WH$ and it is therefore fixed by $K$. 

The following Proposition generalises a result in \cite{broer} which corresponds to the case $K=1$, $\WH=W$.

\begin{proposition}\label{prop:dim1}Assume $\HH$ is finite. If $\dim \Lc=1$, then $X(\HH,K,L)$ is normal if and only if $\Gamma_{\HH,K,L}$ acts non-trivially on $\Lc$.
\end{proposition}
\pf  The restriction map from Proposition \ref{prop:richardson} preserves the grading of the polynomial algebras ${\mathbb C}[\Ec]$ and ${\mathbb C}[\Lc]$ and of their invariant subalgebras. If $\Gamma_{\HH,K,L}$ acts trivially on $\Lc$ then ${\mathbb C}[\Lc]^{\Gamma_{\HH,K,L}}={\mathbb C}[\Lc]$ has terms in degree $1$, whereas 
${\mathbb C}[\Ec]^{\WH}$, and, a fortiori, ${\mathbb C}[\Ec]^{K\ltimes \WH}$, have no components in degree $1$. Thus the restriction map is never surjective. 

If $\Gamma_{\HH,K,L}$ acts non-trivially on $\Lc$, then $V_{\HH,\,{\mathbb C}}$ is necessarily trivial. Also, $\Gamma_{\HH,K,L}$ acts on $V(\Lc)$ by orthogonal transformations preserving the fixed point $p$, so it must act on $V(\Lc)\simeq{\mathbb C}$ as $-1$.
Hence, ${\mathbb C}[\Lc]^{\Gamma_{\HH,K,L}}\simeq{\mathbb C}[V(L)_{\mathbb C}]^{\Gamma_{\HH,K,L}}\simeq {\mathbb C}[t^2]$. The inner product on $V_{\mathbb C}$ is a non-trivial $\WK$-invariant $2$-form, with non-trivial restriction to $V(L)_{\mathbb C}$, so the map is surjective.
\epf

\subsection{Normality in codimension 1}

Here $\HH$ is again arbitrary, admissible. We recall that a variety  is normal in codimension $1$ (unibranch in codimension $1$, respectively) if its non-normality locus (non-unibranchedness locus, respectively) has codimension greater than $1$. In this section we will provide a necessary and sufficient condition for normality  of  $X=X(\HH,K,L)$ in codimension $1$. 

Observe that if $L'=L\cap H$ is a hyperplane in $\HH^L$, then $\WK U_{L'}/\WK$ has codimension $1$ in $X$. Hence, if $X$ is normal, respectively unibranch in codimension $1$, then Proposition \ref{prop:constant} implies that $X$ is normal, respectively unibranch at all points in $\WK U_{L'}/\WK$, whence at all points in $\WK\interiorcp/\WK$. 
Conversely, if $X$ is normal, respectively unibranch, at some point $l'\in U_{L'}$ for every flat $L'=L\cap H$, with $H\in\HH$ and $H\not\supset L$, then
it is normal, respectively unibranch, in codimension $1$.

Let $\Gamma_{\HH,K,L}$ be as in \eqref{eq:gammaK},  let  $\iota_L\colon\Gamma_{\HH,K,L}\to{\rm Aut}(L)$ be induced by restriction to $L$ and let $\WW_{\HH^L}$  be the subgroup    of ${\rm Aut}(L)$ generated by the  reflections with respect to the hyperplanes in $\HH^L$. \label{iotaL}

\begin{lemma}\label{lem:reflection}Let $X=X(\HH,K,L)$ be a stratum and let $l\in U_{L'}\cap\Abar$ for some hyperplane $L'$ in $\HH^L$. Then the finite counterpart $X(\HH_l,K^l,L)$ of $X$ at $\overline{l}$ is normal if and only if the reflection in $L$ with respect to the hyperplane $L'$ lies in $\iota_L(\Gamma_{\HH,K,L})$.  
\end{lemma}
\pf %We consider the finite counterpart $X(\HH_l,K^l,L)$ of $X(\HH,K,L)$ at $l$. 
%Observe that $\Lc'=\bigcap_{H\in\HH_l}H_{\mathbb C}$ and that $\WKl={\WK}_{L'}$ acts trivially on $\Lc'$, so $\WKl$ satisfies the requirements from Section \ref{subsec:basic} for the arrangement $\HH_l$. 
By Remark \ref{rem:reductive} applied to the finite counterpart $X(\HH_l,K^l,L)$ of $X$ at $l$,
\begin{align*}
 \left(X(\HH_l,K^l,L),\,\overline{l}\right)\sim_{se}\left(X(\HH_l',K^l,L/V(L')),\,\overline{l}\right)
\end{align*}
where $\dim(L/V(L'))=1$. 
%we have $Y_l\simeq \Lc'\times \WKl(\Lc/V(\Lc'))/\WKl$, with $\dim(\Lc/V(\Lc'))=1$. 
Since $\HH_l$ is finite, by Proposition \ref{prop:dim1} we have normality at $\overline{l}$, whence everywhere by Remark \ref{rem:real}, if and only if ${\rm Stab}_{\WKl}(\Lc/V(\Lc'))$ acts non-trivially on $\Lc/V(\Lc')$. Any $w\in {\rm Stab}_{\WKl}(\Lc/V(\Lc'))$ fixes $\Lc'$ pointwise because $l\in U_{L'}$.
%$\WKl=\WK^{L'}$. 
Also, ${\rm Stab}_{\WKl}(\Lc/V(\Lc'))\leq \Gamma_{\HH,K,L}$ so $\iota_{L}(w)$ is a non-trivial Euclidean transformation fixing $L'$ pointwise, i.e., the reflection with respect to $L'$. Thus the action is non-trivial if and only if the reflection with respect to $L'$ lies in $\iota_L(\Gamma_{\HH,K,L})$.
\epf

\begin{proposition}\label{prop:nec_codim1}The stratum $X=X(\HH,K,L)$ is normal in codimension $1$ if and only if the following two conditions hold:
\begin{enumerate}[label=(\roman*)]
 \item\label{item:prima} $X$ is unibranch in codimension $1$;
 \item\label{item:seconda} $\WW_{\HH^L}\leq\iota_L(\Gamma_{\HH,K,L})$.
\end{enumerate}
\end{proposition}
\pf We choose $L$ to lie over $\Abar$. If $X$ is normal in codimension $1$, then  \ref{item:prima} holds. Observe that $\Abar\cap \interior$ lies in ${\mathbb C}(\HH^L)$ and that $\WW_{\HH^L}$ is generated by the reflections with respect to the walls of $\AA\cap L$, i.e., the hyperplanes $L'$ in $\HH^L$ lying over $\Abar$. Hence, in order to prove \ref{item:seconda} it is enough to show that  $\iota_L(\Gamma_{\HH,K,L})$ contains the reflections with respect to all  such $L'$. Let $L'$ be one of these and let $l\in U_{L'}\cap\Abar$. By assumption $X$ is normal at $l$. By Corollary \ref{cor:X'}, the stratum $X(\HH_l,K^l,L)$ is normal at $\overline{l}$. By Lemma \ref{lem:reflection}, the reflection with respect to $L'$ lies in $\iota_L(\Gamma_{\HH,K,L})$.  

Conversely, assume conditions \ref{item:prima} and \ref{item:seconda} hold. By Proposition \ref{prop:constant} it is enough to show that $X$ is normal at a point $l$ in $U_{L'}$ for each $L'$ of codimension $1$ in $L$. Since $\Gamma_{\HH,K,L}$ stabilises $\HH^L$, condition (ii) ensure that $\WW_{\HH^L}$ is admissible, so $\Gamma_{\HH,K,L}$ acts transitively on ${\mathcal C}(\HH^L)$ preserving faces. Thus we can replace $L'$ by a $\Gamma_{\HH,K,L}$-translate lying over $\Abar$ and take $l\in U_{L'}\cap\Abar$. By \ref{item:prima} and Corollary \ref{cor:X'}, normality at $\overline{l}$ is equivalent to normality of $X(\HH_l,K^l,L)$, which in turn follows from Lemma \ref{lem:reflection} and \ref{item:seconda}. 
\epf

\begin{remark}\label{rem:semidir}
Assume that condition \ref{item:seconda} from Proposition \ref{prop:nec_codim1} holds.  
\begin{enumerate}[label=(\roman*)]
\item
Since $\Gamma_{\HH,K,L}$ stabilizes $L$ and $\HH$, it stabilizes $\HH^L$ and so does $\WW_{\HH^L}$, hence $\HH^L$ is admissible and $\WW_{\HH^L}$ acts transitively on ${\mathcal C}(\HH^L)$.
\item
If $L$ lies over $\Abar$, then $\Abar\cap \interior$ lies in ${\mathcal C}(\HH^L)$ and $\Abar_L=\Abar\cap L$ is a fundamental domain for the $\WW_{\HH^L}$-action. A standard argument shows that $\iota_L(\Gamma_{\HH,K,L})={\rm Stab}_{\iota_L(\Gamma_{\HH,K,L})}(\Abar_L)\ltimes \WW_{\HH^L}$ and ${\rm Stab}_{\iota_L(\Gamma_{\HH,K,L})}(\Abar_L)=\iota_L({\rm Stab}_K(\Abar_L))$. If $K=1$ this gives $\iota_L(\Gamma_{\HH,L})=\WW_{\HH^L}$.
\end{enumerate}
\end{remark}

\section{The posets ${\mathcal P}(\Sigma(\HH,L))$  and ${\mathcal P}(\Sigma(\HH,K,L))$}\label{sec:sette}

In this section we associate to each stratum $X(\HH,K,L)\subseteq \Ec/\WK$ some subposets of ${\mathcal P}(\HH)$ whose combinatorial properties encode geometric properties of $X(\HH,K,L)$ such as being normal in codimension $1$ and being unibranch. 

\begin{lemma}\label{lem:unions}
Let $L$ be a flat for $\HH$. The subsets $\WK L\cap\Abar$ and $\WH L\cap\Abar$ are unions of closures of maximal faces of dimension $\dim L$. 
\end{lemma}
\pf By construction the two sets  are unions of faces of $\Abar$. The proof for $\WK$ will suffice. If $F=wL\cap\Abar$ for $w\in\WK$, then $F\subset\overline{C}$ for some $C\in{\mathcal C}(\HH^{wL})$. 
For some $\sigma\in\WH$ we have $\sigma C\subset \Abar$ and $\sigma$ fixes $F$ pointwise by Remark \ref{rem:fundamental}. Hence, $F\subset \sigma\overline C\cap\Abar\subset \sigma w L\cap\Abar$ and $\dim |\sigma C|=\dim\sigma wL=\dim L$. \epf

We consider the sets $\Sigma(\HH,K,L)$ and $\Sigma(\HH,L)$ of maximal faces contained in $\WK L\cap\Abar$ and $\WH L\cap\Abar$, respectively. It follows from the proof of Lemma \ref{lem:unions} that \label{Sigma}
\begin{align*}\Sigma(\HH,K,L)&:=\{w\interior\cap\Abar~:~w\in\WK,\;\;wL \mbox{ lies over }\Abar \},\\
\Sigma(\HH,L)&:=\{w\interior\cap\Abar~:~w\in \WH,\;\;wL \mbox{ lies over }\Abar\}.\end{align*}
%\begin{equation*}\Sigma_K:=\{~|~wL\in \Xi_K,\,{\mathcal}\},\quad \Sigma_\HH=\{wL\cap\Abar~|~wL\in \Xi_\HH\}\end{equation*}
By construction $\Sigma(\HH,K,L)$ is $K$-stable and equals $K\Sigma(\HH,L)$. The set $\Sigma(\HH,K,L)$ uniquely determines $X(\HH,K,L)$, so the collection of sets of this form parametrises strata in $\Ec/\WK$. 
The proof of Lemma \ref{lem:unions} shows that any real point in $X(\HH,K,L)$ is represented by a point in the closure of some face in $\Sigma(\HH,K,L)$. In particular, by Remark \ref{rem:real} minimal strata are represented by points in minimal dimensional faces in the closure of a face in $\Sigma(\HH,K,L)$. For this reason, it is important to consider the induced subposets of ${\mathcal P}(\HH)$, consisting of faces contained in $\WK L\cap\Abar$ and $\WH L\cap\Abar$, respectively.\label{posets}
\begin{align*}
 &{\mathcal P}(\Sigma(\HH,K,L)):=\{F\in{\mathcal P}(\HH^L)~:~F\subset \overline{F'}, \text{ for some } F'\in\Sigma(\HH,K,L) \},\\
&{\mathcal P}(\Sigma(\HH,L)):=\{F\in{\mathcal P}(\HH^L)~:~F\subset\overline{F'} \text{ for some } F'\in\Sigma(\HH,L) \}.\end{align*}

In order  to verify normality and unibranchedness of $X(\HH,K,L)$ it is enough to verify it at a point in all  minimal dimensional faces in ${\mathcal P}(\Sigma(\HH,K,L))$. If $K$ acts transitively on $\Sigma(\HH,K,L)$ it is enough to verify unibranchedness and normality at the minimal dimensional faces contained in the closure of one face in $\Sigma(\HH,K,L)$.

\begin{lemma}\label{lem:fsigma}Let $L$ be a flat.
There is a dimension-preserving surjective poset map $f_\Sigma\colon{\mathcal P}(\HH^L)\to{\mathcal P}(\Sigma(\HH,L))$ induced from a piecewise-linear surjective map $L\to \Abar$.
%\begin{enumerate}
%\item\label{item:fsigma} 
%\item\label{item:equiv-real} There is a well-defined surjective simplicial map $\WW\Sigma_{v+\zz}\to\Sigma_J$ which is $\WW$-invariant and induces a bijection of sets
%$\WW\Sigma_{v+\zz}/\WW\to\Sigma_J$.
%\end{enumerate}
\end{lemma}
\pf We consider the piecewise-linear map $f\colon L\to \Abar$ associating to each $l\in L$ the unique $p\in\Abar\cap \WH l$. 
All points lying in the same face $F$ in ${\mathcal P}(\HH^L)$ are mapped to points in a face contained in $\Abar$ of the same dimension as $F$. In particular, any chamber of $L$ is mapped to a unique face of $\Sigma(\HH,L)$. By construction, $f$ preserves inclusion of closures, inducing the sought poset map $f_\Sigma$. If $F$ is a face in $\Sigma(\HH,L)$, then  $|F|=wL$ for some $w\in \WH$ and $wL$ lies over $\Abar$. Moreover, $w^{-1}F\in\HH^L$ so the points in $F$ lie in the image of $f$ whence $f_\Sigma$ is surjective. \epf

\begin{corollary}\label{cor:connected}
For any two distinct maximal faces $F$, $F'$ in ${\mathcal P}(\Sigma(\HH,L))$ there exists a gallery of maximal faces 
beginning at $F$ and ending at $F'$.
\end{corollary}
\pf Let $F,\,F'$ be maximal faces in  ${\mathcal P}(\Sigma(\HH,L))$. Since $f_\Sigma$ is surjective and dimension preserving, there exist maximal faces $C,\,C'$ in ${\mathcal P}(\HH^L)$ such that $f_\Sigma(C)=F$, $f_\Sigma(C')=F'$. Since $L\setminus \bigcup_{\genfrac{}{}{0pt}{}{F\in {\mathcal P}(\HH^L)}{ \text{codim}_L F\geq2}}F$ is path connected, any path from a point in $C$ to one in $C'$ determines a gallery of maximal faces $C_0=C,\,\ldots,\, C_r=C'$ in ${\mathcal P}(\HH^L)$. The required gallery is obtained applying $f_\Sigma$ to this sequence and removing possible repetitions occurring for those $i$ such that  $f_\Sigma(C_i)=f_\Sigma(C_{i+1})$. 
\epf

\begin{remark}\label{rem:stationary}\begin{enumerate}[label=(\roman*)]
\item Let $C_1\neq C_{2}$ be adjacent maximal faces in ${\mathcal P}(\HH^L)$. Assume $f_\Sigma(C_1)=f_\Sigma(C_{2})$ 
and let $w_1,\,w_2\in \WH$ such that $w_1C_1=w_2C_{2}=f_\Sigma(C_1)\subset\Abar$. Then, $w_2^{-1}w_1\in\Gamma_{\HH,L}$ and 
$w_1l=w_2l$ for every $l\in \overline{C_1}\cap\overline{C_{2}}$. Therefore, the Euclidean map $\iota_L(w_2^{-1}w_1)$ is necessarily the reflection with respect to the wall in $\HH^L$ separating $C_1$ and $C_2$. 
\item Let $C_i$ for $i=0,\,\ldots,\, m$ form a gallery ${\mathcal G}$ of maximal faces in ${\mathcal P}(\HH^L)$ and let $w_i\in \WH$ be such that $w_iC_i=F_i\subset\Abar$ for every $i$. There is always a gallery ${\mathcal G}'$ of maximal faces $C'_i$  in ${\mathcal P}(\HH^L)$ for $i=0,\,\ldots m'$ with $m'\leq m$ satisfying $f_\Sigma(C_0)=f_\Sigma(C_0')$, $f_\Sigma(C_m)=f_\Sigma(C'_{m'})$ and $f_\Sigma(C'_i)\neq f_\Sigma(C'_{i+1})$ for every $i$. Indeed,
if  $j$ is the minimum index for which  $f_\Sigma(C_j)=f_\Sigma(C_{j+1})$, we can replace ${\mathcal G}$ by the shorter gallery $C_0,\,\ldots,\, C_j=w_{j+1}^{-1}w_j C_{j+1},\, w_{j+1}^{-1}w_j C_{j+2}, \,\ldots, w_{j+1}^{-1}w_j C_{m}$.
Then $f_\Sigma(C_0)=F_0$ and $f_\Sigma(w_{j+1}^{-1}w_j C_{m})=F_m$. Iterating this procedure gives ${\mathcal G}'$.
\end{enumerate}
\end{remark}

\begin{lemma}\label{lem:condition}Let $L$ lie over $\Abar$. If $\WW_{\HH^L}\leq\iota_L(\Gamma_{\HH,K,L})$, then for any two faces $F_1,\,F_2$ in $\Sigma(\HH,L)$ there is $k\in K$ such that $kF_1=F_2$. In particular, if $X(\HH,K,L)$ is normal in codimension $1$, then $K$ acts transitively on $\Sigma(\HH,K,L)$.
\end{lemma}
\pf  Assume $\WW_{\HH^L}\leq\iota_L(\Gamma_{\HH,K,L})$ and let $F_1,\,F_2$ in $\Sigma(\HH,L)$. For $i=1,2$ let $C_i\in{\mathcal P}(\HH^L)$  be such that $f_\Sigma(C_i)=F_i$ and let $w_i\in\WH$ be such that $w_iC_i=F_i$. By Remark \ref{rem:semidir} the group $\iota_L(\Gamma_{\HH,K,L})$ acts transitively on the set of maximal faces in ${\mathcal P}(\HH^L)$, so there is $kw\in K\ltimes\WH$ such that $kw C_1=C_2$. Then, $kw w_1^{-1}F_1=w_2^{-1}F_2$ and $k(k^{-1}w_2kww_1^{-1})F_1=F_2\in\Abar$. Remark \ref{rem:fundamental} applies and $kF_1=F_2$. Last statement follows from Proposition \ref{prop:nec_codim1} because $\Sigma(\HH,K,L)=K\Sigma(\HH,L)$.\epf

%The following Lemma is a partial converse of Lemma \ref{lem:condition}.
%\begin{lemma}\label{lem:partial}
%Let $L$ lie over $\Abar$ with $F=\interior\cap\Abar$ and let $k\in K$. Assume that $F$ and $kF\in\Sigma_K$ are adjacent and  
%the face $F''$ separating them is fixed pointwise by $k$. Then, $\iota(\Gamma_K)$ contains the reflection with respect to $|F''|$.  
%\end{lemma}
%\pf Let $F'$ be the chamber of $\HH^L$ separated from $F$ by the wall $H=|F''|$. Then $f_\Sigma(F')$ contains $\overline{F}\cap\overline{F'}\subset \Abar$, so either $f_\sigma(F')=F$ or 
%$f_\sigma(F')=kF$. In the first case we apply Remark \ref{rem:stationary} (i). In the second case, we consider $w\in\WH$ such that $wF'=kF$. Then, $k^{-1}w\in\Gamma_K$ and its action on $L$ is orthogonal and fixes $|\overline{F}\cap\overline{F'}
%|=H$ pointwise by Remark \ref{rem:fundamental}.
%\epf

The following Lemma shows how to describe $\Omega(\HH,K,L)_{l}$ for $l\in\Abar$ in terms of faces in $\Sigma(\HH,K,L)$. For $l\in L\cap\Abar$ we set \label{Fl}
\begin{align*}
{\mathcal F}_l=\{F\in\Sigma(\HH,K,L)~:~l\in\overline{F}\}
\end{align*}

\begin{lemma}\label{lem:set}For $L$ a flat lying over $\Abar$ and $l\in L\cap\Abar$
%and let ${\mathcal F}_l=\{F\in\Sigma(\HH,K,L)~:~l\in\overline{F}\}$. Then
there holds
\begin{align*}\Omega(\HH,K,L)_{l}=\{w |F|_{\mathbb C}~:~w\in\WHl,\,F\in{\mathcal F}_l\}.\end{align*}
\end{lemma}
\pf By construction we have the inclusion $\supseteq$. We prove $\subseteq$. Let $w\in\WK$ such that $l\in wL\subset w\Lc$ and let $C\in{\mathcal C}(\HH^{wL})$ with $l\in\overline{C}$. By \cite[V.3.3, Remarque 1]{bourbaki} there is $\sigma\in\WH$ such that $\sigma C\subset\Abar$. Since $l,\,\sigma l\in\Abar$, we have $\sigma\in\WHl$. In addition, 
\begin{align*}l\in\sigma wL=\sigma|C|=|\sigma C|\subseteq |(\sigma wL)\cap\Abar|\subseteq \sigma wL.\end{align*} Hence $\sigma w L$ lies over $\Abar$ so $F':=\Abar\cap \sigma w\interior$ lies in ${\mathcal F}_l$ and $w\Lc=\sigma^{-1}|F'|_{\mathbb C}$. 
\epf

We aim at giving a characterisation of unibranchedness and normality in codimension $1$ in terms of the $K$-action on ${\mathcal P}(\Sigma(\HH,K,L))$. 

\begin{lemma}\label{lem:id}
Let $L'\subset L$ be flats lying over $\Abar$, with $F=\interior\cap\Abar$ and $F'=\interiorp\cap\Abar$. If
% and assume that the following condition holds:
\begin{equation}\label{eq:condition}
%\{kF~:~k\in K^{L'}\}=
\{F''\in\Sigma(\HH,K,L)~:~F'\subset\overline{F''}\} \textrm{ is a single $K^{L'}$-orbit}, 
\end{equation}
then, $X(\HH,K,L)$ is unibranch at all points in $\interiorcp$. 
\end{lemma}
\pf By Proposition \ref{prop:constant} it is enough to show that if \eqref{eq:condition} holds, then $X(\HH,K,L)$ is unibranch at some $l\in F'\cap U_{L'}$. Since $L'$ is the minimal flat containing $l$, if $F''\in{\mathcal F}_l$
%
%
%$l\in \overline{F''}$ for some $F''\in\Sigma(\HH,K,L)$, 
then $\overline{F''}\supset \interiorp\cap\Abar=F'$.
Hence Lemma \ref{lem:set} and \eqref{eq:condition} give 
\begin{equation*}\Omega(\HH,K,L)_{l}=\{w k|F|_{\mathbb C}~:~w\in\WHl,\,k\in K^l\}\end{equation*}
where we have used that $K^l=K^{L'}$. We conclude by Corollary \ref{cor:condition_unibranch}.
\epf

\begin{lemma}\label{lem:idd}
Let $L'\subset L$ be flats lying over $\Abar$. Assume $\WW_{\HH^L}\leq\iota_L(\Gamma_{\HH,K,L})$. Then, $X(\HH,K,L)$ is unibranch at $\interiorcp$ if and only if condition \eqref{eq:condition} holds for $F'=\interiorp\cap\Abar$.
\end{lemma}
\pf We prove the converse of Lemma \ref{lem:id}. 
Assume $X(\HH,K,L)$ is unibranch at $l\in U_{L'}\cap\Abar$ and
%for $L'\subset L$ lying over $\Abar$. 
let $F''\in \Sigma(\HH,K,L)$ be such that $F'\subset \overline{F''}$. Then $L'':=|F''|=\sigma L$ for some $\sigma=w^{-1}k^{-1}\in \WK$. By Corollary \ref{cor:condition_unibranch} we may take $k\in K^l=K^{L'}$ and $w\in \WHl=\WW_{\HH,K}^{L'}$. Now $kwF''$ and $F$ lie in ${\mathcal C}(\HH^L)$ so by Remark \ref{rem:semidir} (i), there exists $k_1w_1\in\iota_L^{-1}(\WW_{\HH^{L}})$ such that $k_1w_1kwF''=F$.

Remark \ref{rem:fundamental} applied to $\Abar$ gives $k_1w_1kwF''=k_1kF''=F$ and  
$k_1w_1x=k_1x$, and $k_1kx=k_1w_1kwx=k_1w_1x$ for any $x\in F'$. The same argument applied to the fundamental domain $\overline{F}$ for the action of $\WW_{\HH^L}$ gives $k_1w_1x=x$ for $x\in F'\subset \overline{F}$. Hence, $k_1kx=k_1w_1x=x$ and $k_1k \in K^{L'}$ giving \eqref{eq:condition}.
\epf

\begin{proposition}\label{prop:equivalent}Let $L$ be a flat lying over $\Abar$ and let $F=\interior\cap\Abar$. Then $X(\HH,K,L)$ is normal in codimension $1$ if and only if \eqref{eq:condition} holds for every $L'\subset L$ of codimension $1$ lying over $\Abar$. 
\end{proposition}
\pf If $X(\HH,K,L)$ is normal in codimension $1$, then Proposition \ref{prop:nec_codim1} and  Lemma \ref{lem:idd} give $\WW_{\HH^L}\leq\iota_L(\Gamma_{\HH,K,L})$ and \eqref{eq:condition} holds for every $L'\in\HH^L$ lying over $\Abar$.
%by Proposition \ref{prop:nec_codim1} and  Lemma \ref{lem:idd}.
%Let $F_1\in \Sigma_K$ such that $F'\subset \overline{F_1}$. Then $|F_1|=L_1=\sigma L$ for some $\sigma=w^{-1}k^{-1}\in \WK$. By Corollary \ref{cor:condition_unibranch} we may take $k\in K_l=K_{L'}$ and $w\in \WHl=\WW_{L'}$. Then $kwF_1$ is a chamber of $L$, so there exists $k_1w_1\in\iota^{-1}(\WW_{\HH^{L}})$ such that $k_1w_1kwF_1=F$. By Remark \ref{rem:fundamental} applied to $\Abar$ we have $kwx=kx=x$ for every $x\in F'$ and $k_1w_1kwx=k_1kx$. The same argument applied to the fundamental domain $\overline{F}$ for $\WW_{\HH^L}$ gives $k_1w_1x=x$. Hence, $k_1k \in K_{L'}$. 

Conversely, assume that  \eqref{eq:condition} holds for every $L'\in \HH^L$ lying over $\Abar$. We prove condition \ref{item:seconda} of Proposition \ref{prop:nec_codim1}. It is enough to prove that all reflections with respect to a wall of the chamber $F$ in $\HH^L$ lie in the image of $\iota_L$. Let $L'=H\cap L$ be such a wall, let $C\in{\mathcal C}(\HH^L)$ be the chamber adjacent to $F$ on the other side of $L'$, and let $F'=\interiorp\cap\Abar$. Then $F_1:=f_\Sigma(C)$ is a chamber lying in $\Sigma(\HH,L)\subset\Sigma(\HH,K,L)$ containing $F'$ in its closure, so  $F_1=wC$ for some $w\in \WH$ with $w\in\WW_\HH^{L'}$ by Remark \ref{rem:fundamental}. If $F_1=F$, then the reflection with respect to the wall $H\cap L$ lies in $\iota_L(\Gamma_{\HH,L})$ by Remark \ref{rem:stationary}. If $F_1\neq F$, there is $k\in K^{L'}$ such that $kw C=F$. Then $kw\in \iota_L(\Gamma_{\HH,K,L})$, it fixes $L'$ pointwise, so it is the sought reflection.  

We prove condition \ref{item:prima} of Proposition \ref{prop:nec_codim1}. For those $L'\in\HH^L$ lying over $\Abar$ this follows from Lemma \ref{lem:id}. 
For those $L'\in\HH^L$ not lying over $\Abar$, it is enough to prove unibranchedness at a real element $l\in U_{L'}$. Let $C\in{\mathcal P}(\HH^L)$ with $l\in C$.  Then, there is $w\in\iota_L(\Gamma_{\HH,K,L})$ such that $wl\in wC\subset\Abar$, so $wL'$ lies over $\Abar$ and $X(\HH,K,L)$ is normal at the class $\overline{l}=\overline{wl}$.
\epf

\section{Coxeter classes}\label{sec:combinatorial}

In this section we show how to compute $\Sigma(\HH,L)$ in terms of subsets of the set $S$ of Coxeter generators of $\WH$ given by the reflections with respect to the walls of some component $\AA_{(j)}$ of $\AA$. We identify $S$ with the set of nodes $\{N_1,\,\ldots,\,N_n\}$ of the Coxeter graph of $\WH$.  \label{S}

Let $F\in{\mathcal P}(\HH)$ such that $F\subset\Abar$. We associate to $F$ the subset $S_F$ of $S$ consisting of the reflections with respect to the  walls of $\AA$ containing $F$. The corresponding set of walls is denoted by $M_F$. \label{emmeffe} By construction, $|F|=\bigcap_{H\in M_F}H$. We also set $\hat{S}_F:=S\setminus S_F$. The parabolic subgroup $\WW_{S_F}\leq\WH$ generated by the reflections in $S_F$ is $\WW_{\HH}^{|F|}$, the pointwise stabiliser of $|F|$. 

If $\WH$ is not irreducible, then $F=\prod_{i=1}^rF_{(i)}$ is a product of faces corresponding to each component of $E$ and the subset $S_F$ is compatible with this decomposition. The sets $S_F$ run through all subsets of  $S$ that do not contain a whole affine component, i.e., subsets for which $\WW_{S_F}$ is  {\em finite}. 

Since $K$ acts on the faces of ${\mathcal P}(\HH)$ contained in $\Abar$, it acts on the collection of walls of the form $M_F$ for some $F\in \Sigma(\HH,K,L)$ and thus on the corresponding collection of subsets of $S$ of the form $S_F$.
%and thus on ${\mathcal S}_{\HH,K,X}$
%To any stratum $X(\HH,K,L)$, we associate the set: \label{SX}
%\begin{align*}{\mathcal S}_{\HH,K,L}:=\{J\subset S~:~J=J_F \textrm{ for some }F\in\Sigma(\HH,K,L)\}\end{align*}
%Since $K$ acts on the faces of ${\mathcal P}(\HH)$ contained in $\Abar$, it acts on the collection of  walls $M_F$ and thus on ${\mathcal S}_{\HH,K,X}$. 

Observe that $K$ does not necessarily preserve components of $E$ and $S$.  Also, if $1\neq w\in\WH$, then $w$ does not preserve the set of walls of $\AA$. However, if for some wall $H$ the hyperplane $wH$ is again a wall of $\AA$, then we say that the image of the corresponding node through $w$ lies in $S$ and it is the node associated to $wH$. Following \cite{GP}, for a subset  $S_F\subset S$ we define  the {\em Coxeter class of $S'$} as the subset \label{CoxeterJ}
\begin{align*}[S_F]:=\{S'\subset S~:~wS_F=S' \textrm{ for some } w\in\WH\}.\end{align*}
\begin{proposition}\label{prop:Gamma'}
Let $L$ be a flat lying over $\Abar$ with $\AA_L=\interior\cap\Abar$ in $\Sigma(\HH,L)$.  The assignment $F\mapsto S_{F}$ sets a bijection between  $\Sigma(\HH,L)$ and  the Coxeter class $[S_{\AA_L}]$. \end{proposition}
\pf If $F'\in\Sigma(\HH,L)$, then $L':=|F'|=wL$ for some $w\in\WH$. Conjugation by $w$ maps $\WW_{\HH}^{L}=\WW_{S_{\AA_L}}$ to $\WW_{\HH}^{L'}=\WW_{S_{F'}}$. We claim that there is $\sigma\in\WH$ mapping $S_{\AA_L}$ to $S_{F'}$. For $\WH$ finite this is  \cite[Corollary 2.1.13]{GP}. In the same spirit, if $\WH$ is infinite, conjugation by the minimal length representative $\sigma\in\WW_{S_{F'}}w\WW_{S_{\AA_L}}$ maps length $1$ elements of $\WW_{S_{\AA_L}}$ to length $1$ elements of $\WW_{S_{F'}}$, \cite[Chap. IV, Exercise \S 1 n. 3]{bourbaki}. Hence, the assignment $F'\mapsto S_{F'}=\sigma S_{\AA_L}$ has image in  $[S_{\AA_L}]$. By construction, the map is injective.  Let  $\sigma\in \WH$ with $\sigma S_{\AA_L}\subset S$. Then, for $L':=|\sigma S_{\AA_L}|=\sigma L$ and $F':=\interiorp \cap\Abar$ in $\Sigma(\HH,L)$, we have $\sigma S_{\AA_L}=S_{F'}$ giving surjectivity.\hfill$\Box$

\begin{corollary}\label{cor:parametrization}Strata in $\Ec/\WH$ are in bijection with Coxeter classes of subsets $S'\subset S$ generating \emph{finite} parabolic subgroups of $\WH$. Strata in $\Ec/\WK$ are in bijection with sets of the form $KZ$ where $Z$ is a Coxeter class as above.  
\end{corollary}

We produce an algorithm to compute $\Sigma(\HH,L)$, and therefore $\Sigma(\HH,K,L)=K\Sigma(\HH,L)$  starting from an element therein.

\begin{proposition}\label{prop:algo}Let $L$ be a stratum lying over $\Abar$ and let $F\in\Sigma(\HH,L)$. A face $F'$ lies in $\Sigma(\HH,L)$ if and only if there exists a gallery $\{F_i\in \Sigma(\HH,L),\, 0\leq i\leq m\}$ with
$F_0=F$, $F_m=F'$ and such that $S_{F_i}$ is mapped to $S_{F_{i+1}}$ by the longest element in the finite group $\langle \WW_{S_{F_i}},\WW_{S_{F_{i+1}}}\rangle\cap{\WH}_{(j)}$, where $j$ is the unique index on which the components ${F_i}_{(j)}$ and ${F_{i+1}}_{(j)}$ differ. 
\end{proposition}
\pf If  such a gallery exists, then $F'$ is a face in $\Sigma(\HH,L)$ by Proposition \ref{prop:Gamma'}.

Let $F'\in \Sigma(\HH,L)$ and let $C,\,C'$ be faces in ${\mathcal P}(\HH^L)$ such that $f_\Sigma(C)=F$, $f_\Sigma(C')=F'$. By Remark \ref{rem:stationary} (ii), there is a gallery $C_0=C,\,\ldots,\,C_m=C'$  of maximal faces in ${\mathcal P}(\HH^L)$ such that $F_i:=f_\Sigma(C_i)\neq f_\Sigma(C_{i+1})=:F_{i+1}$. Let $H_i=\bigcap_{H_\in M_{F_i}\cup M_{F_{i+1}}}H$ be the unique wall separating $F_i$ and $F_{i+1}$, let $L_i=|F_i|$ and let $w_i\in \WH$ be such that $w_iC_i=F_i$, so $w_iL=L_i$. Then, for $\sigma_i:=w_{i+1}w_i^{-1}$ we have $\sigma_i L_i=L_{i+1}$ and, since $\sigma_i(H_i\cap\Abar)\subset \Abar$, the element $\sigma_i$ acts as the identity on $H_i$. Hence, $\sigma_i$ lies in the group generated by the reflections with respect to all hyperplanes containing $H_i$, i.e., $\langle \WW_{S_{F_i}},\WW_{S_{F_{i+1}}}\rangle$. By construction, $\sigma_i(S_{F_i})=S_{F_{i+1}}$. Our procedure shows that $w_i$ and $w_{i+1}$ can be chosen in the component ${\WH}_{(j)}$, hence so does $\sigma_i$.
 We claim that the parabolic subgroup $\langle \WW_{S_{F_i}},\WW_{S_{F_{i+1}}}\rangle\cap {\WH}_{(j)}$ of $\WH$ is finite.
 Indeed, it could be infinite only if the Coxeter graph of ${\WH}_{(j)}$ were the underlying graph of an extended Dynkin diagram of a simple Lie algebra, and $\langle \WW_{S_{F_i}},\WW_{S_{F_{i+1}}}\rangle\cap {\WH}_{(j)}={\WH}_{(j)}$. If this were the case, the $j$-th component ${F_i}_{(j)}$ of the face ${F_i}$ would be a point, and likewise for $F_{(j)}$ and for the component of $L$ in $E_{(j)}$. 
By definition of $\Sigma(\HH,L)$, we would have ${F_i}_{(j)}=F_{(j)}={F_{i+1}}_{(j)}$, a contradiction. 
We apply thus \cite[Proposition 2.3.2 (i)]{GP} to $\langle \WW_{S_{F_i}},\WW_{S_{F_{i+1}}}\rangle\cap {\WH}_{(j)}$ to see that  $\sigma_i$ can be chosen to be the longest element therein. 
\hfill$\Box$

\begin{remark}
 If $\WH$ is finite, Propositions \ref{prop:Gamma'} and \ref{prop:algo} give \cite[Theorem 2.3.3]{GP} in force of \cite[Proposition 2.3.2 (i)]{GP}. 
 \end{remark}

Proposition \ref{prop:algo} indicates a procedure to compute $\Sigma(\HH,L)$.
\begin{corollary}\label{cor:algo}
Algorithm F in \cite[Section 2.3]{GP} can be applied to compute $\Sigma(\HH,L)$ from $F\in\Sigma(\HH,L)$ also when $\WH$ has affine components.
\end{corollary}
\pf It is enough to verify this for $\WH$ irreducible. If $F$ is a point, then $\Sigma(\HH,L)=\{F\}$. If $F$ is not a point, we use Proposition \ref{prop:algo} instead of  \cite[Theorem 2.3.3]{GP}. For every node $N_j$ in $\hat{S}_F$ we apply  to $S_F$ the longest element $w_{0,j}$ of the parabolic subgroup of $\WH$ generated by the reflections corresponding to $N_j$ and  $S_F$. If $S_F\neq w_{0,i}S_F\subset S$, then $w_{0,j}S_F=S_{F_j}$ for some $F_j\in \Sigma(\HH,L)$ and we add it to $\Sigma(\HH,L)$, otherwise we do not add faces to $\Sigma(\HH,L)$. This way we get a new set of elements in $\Sigma(\HH,L)$. Iterating the procedure to all subsets $S'$ constructed this way and all nodes in $\hat{S'}$, we obtain the full set $\Sigma(\HH,L)$.  
\epf

\section{Normality for $K=1$}\label{sec:groupcase}

In this section $K=1$, so $U_L=\interiorc$ for any flat $L$ and the geometry of a stratum is constant on $\WH\interiorc/\WH$. Remark \ref{rem:real} guarantees that  up to smooth equivalence, strata are products of strata corresponding to the irreducible factors of $\WH$ and they are normal, respectively unibranch, respectively smooth if and only if each factor is so. Normality and smoothness when $\WH$ is finite has been dealt with in \cite{richardson,broer, DR}. Thus, it remains to consider the case of $\WH$ effective, irreducible and affine.

The following result was firstly observed in \cite{broer} for $\WH$ finite as a consequence of Chevalley-Shephard-Todd's theorem.  
\begin{corollary}\label{cor:equivalence}
A stratum $X=X(\HH,L)$ is smooth if and only if it is normal. 
\end{corollary}
\pf Smoothness implies normality so we only need to prove the converse. It is enough to prove it for $\WH$ effective and irreducible. If $X$ is normal, then it is normal at all points in minimal, i.e., $0$-dimensional strata by 
Corollary  \ref{cor:dim0}. Each of these is represented by an $l$ in $\Abar$  so  Corollary \ref{cor:X'}  gives $(X,\overline{l})\sim_{se}(X(\HH_l,L),\,\overline{l})$. The right hand side is a normal stratum for a finite reflection group, so it is smooth by \cite[Theorem 3.1]{broer}. Hence, the singular locus of $X$ does not contain $0$-dimensional strata. 
However, in our situation $U_{L'}=\interiorcp$ for every flat $L'\subset L$.  By Proposition \ref{prop:constant} the singular locus, if non-trivial, would contain $0$-dimensional strata because it is closed. Hence, $X$ is smooth.
\epf

\begin{lemma}\label{lem:due-uni}Let $L$ be  flat lying over $\Abar$. If $\#\Sigma(\HH,L)\leq2$, then $X(\HH,L)$ is unibranch.  
\end{lemma}
\pf Lemma \ref{lem:unions} and Proposition \ref{prop:constant} (i) imply that it is enough to verify unibranchedness at points in the closure of faces in $\Sigma(\HH,L)$. 
If $\# \Sigma(\HH,L)=1$,  then \eqref{eq:condition} is trivially satisfied, so the statement follows from  Lemma \ref{lem:id}.
Assume now $\Sigma(\HH,L)=\{F_0, F_1\}$, with $F_0\subset L$ and let $F'$ be the face contained in $\overline{F_0}\cap\overline{F_1}$, whose existence is guaranteed by Corollary \ref{cor:connected}.  We will prove unibranchedness at $l\in\overline{F_0}\cup\overline{F_1}$ using Lemma \ref{lem:set} and Corollary \ref{cor:condition_unibranch}.  If $|{\mathcal F}_l|=1$, then this is immediate. If ${\mathcal F}_l=\{F_0,\,F_1\}$, then $l\in\overline{F'}$.  By Remark \ref{rem:stationary} (ii), there is a gallery $C_i$, for $i=1,2$ of chambers in $L$ such that $C_0=F_0$ and $f_\Sigma(C_1)=F_1$, so $F_1=\sigma C_1$ for some $\sigma\in\WH$ fixing $F'$ pointwise.  By Lemma \ref{lem:set} 
\begin{equation*}\Omega(\HH,1,L)_{l}=\{w|F_0|_{\mathbb C}~:~ w\in \WHl\}\cup \{w\sigma|C_1|_{\mathbb C}~:~ w\in \WHl\}.
\end{equation*} Since $L=|F_0|=|C_1|$, the statement follows from Corollary \ref{cor:condition_unibranch}.\epf

Next Lemma  translates normality in codimension $1$ and unibranchedness into statements concerning $\Sigma(\HH,L)$. 

\begin{lemma}\label{lem:effe}
Let $L$ be a flat for $\HH$. The following statements are equivalent
\begin{enumerate}[label=(\roman*)]
\item\label{item:uno} $\iota_L(\Gamma_{\HH,L})=\WW_{\HH^L}$.
\item\label{item:due} $\Sigma(\HH,L)$ has a unique element.
\item\label{item:quattro} $X(\HH,L)$ is normal in codimension $1$.
\end{enumerate}
\end{lemma}
\pf Condition \ref{item:uno} implies \ref{item:due} by Lemma \ref{lem:condition}.  On the other hand, if  \ref{item:due} holds, then $X(\HH,L)$ is normal in codimension $1$ by Proposition \ref{prop:equivalent}. 
Finally, (iii) implies (i) by Proposition \ref{prop:nec_codim1} and Remark \ref{rem:semidir} (ii). \epf

\begin{remark}When $\WH$ is finite, the property of $X(\HH,L)$ being normal in codimension 1 is equivalent to $S_F$ being self-opposed in the terminology of \cite[2.3.5]{GP}, see also Corollary \ref{cor:algo}. It is also equivalent to the equality of exponents in  \cite[Theorem 2.1]{DR}. \end{remark}

\begin{remark}
If a Levi subgroup $L$ of a parabolic subgroup of $G$ supports a cuspidal local system as in \cite[2.4]{lusztig-inventiones}, the results in \cite[\S 9.2]{lusztig-inventiones} show that the quotient $W\zz(\ll)/W$, for $\ll={\rm Lie}(L)$, is normal in codimension $1$ \cite[\S 9.2]{lusztig-inventiones}. A list of such $L$ is to be found in {\em loc. cit.} or in \cite[\S 2.13]{lusztig-hecke}. Such quotients are also normal, although they do not exhaust the list of normal strata in $\gg/\!/G$.
\end{remark}

\subsection{List of normal Jordan strata in $G_{sc}/\!/G_{sc}$}

In this Section $\WH\simeq W_{aff}$ is irreducible and acts effectively on $E$,  i.e., we are studying strata in $G/\!/G$ for $G=G_{sc}$ simple and simply connected. Here ${\mathcal A}$ is the fundamental alcove and the Coxeter graph of $\WH$ is the underlying graph of the extended Dynkin diagram of $W$. A face $F\subset\Abar$ is the simplex generated by the vertices corresponding to the nodes in $\hat{S}_F$. 

We will provide an answer to the normality problem for Jordan strata in $G/\!/G$ and we collect here how some of the objects we studied are translated into this language. Let $L$ be a stratum  lying over $\Abar$ and let $\overline{J}/\!/G$ be the stratum corresponding to $X(\HH,L)$ through the identification in Proposition \ref{prop:stratifications}. If $F\in \Sigma(\HH,L)$, then $S_F$ corresponds to the (positive and negative)  root subgroups generating with $T$ the centraliser of a representative of the semisimple Jordan class $J$ in $G$. Minimal strata contained in $\overline{J}/\!/G\simeq X(\HH,L)$ are quotients of closures of those Jordan classes contained in $\overline{J}$ and consisting of a unique semisimple class. We recall that the semisimple classes that are themselves Jordan classes in $G$ are precisely those with semisimple connected centraliser, i.e., the isolated semisimple classes, \cite[Definition 2.6]{lusztig-inventiones}. For any $l\in \Abar\cap L$, the finite counterpart $X(\HH_l,L)$ is a stratum for a finite reflection group $\WW_{\HH_l}$.  Through the identification in Proposition \ref{prop:Lie-stratification}, and fixing $l$ as an origin of $\Ec$,  we see that $X(\HH_l,L)$ corresponds to the Jordan stratum $\mathfrak{X}$ for $\mathfrak{c}_{\gg}(e(l))$ given by $W_{e(l)}V(L)_{\mathbb C}/W_{e(l)}$, where $W_{e(l)}$ is the Weyl group of $\mathfrak{c}_{\gg}(e(l))$. If $l$ lies in a minimal stratum, then $e(l)$ is isolated and $\mathfrak{c}_{\gg}(e(l))$ is a semisimple Lie algebra.

We produce first the list of strata that are normal in codimension $1$. 
\begin{proposition}\label{prop:codim1} Let $X=X(\HH,L)$ be a stratum in $\Ec/W_{aff}$ and let $F$ be a maximal face in $\Sigma(\HH,L)$.
Then $X$ is normal in codimension 1 if and only if $S_F=\emptyset$, or $\#S_F=\#{\Delta}$ or it is as follows:
\begin{itemize}
 \item[$A_{n}$ :] of type $dA_{h}$ with $n+1=d(h+1)$, $h\geq 1$, $d\geq 2$;
  \item[$B_n$ :] of type $D_{m_0}+ dA_h + B_{n_0}$ with $n=m_0+n_0+d(h+1)$ and either $m_0\geq 2$, $n_0\geq 0$, $h\geq 0$, or else  $m_0=0$, $n_0\geq 0$, $h=0$ or odd;
  \item[$C_n$ :] of type $C_{m_0}+C_{n_0}+dA_h$ with $m_0,n_0,h\geq0$, $n=m_0+n_0+d(h+1)$;
  \item[$D_n$ :] of type $D_{m_0}+D_{n_0}+dA_h$ with  $n=m_0+n_0+d(h+1)$ and either $m_0,n_0\geq2$ and $h\geq0$, or else $m_0n_0=0$ and $h=0$ or odd;
  \item[$E_6$ :]  of type $A_5$ (there are three such subsets), $D_4$, $4A_1$, $2A_2$ (there are three such subsets);
  \item[$E_7$ :] of type $E_6$, $D_6$ (there are two such subsets), $D_5+ A_1$ (there are two such subsets), $D_4+2A_1$, $2A_3$, $3A_2$, $A_3+3A_1$ (there are two such subsets),  $D_4+A_1$ (there are two such subsets), $5A_1$, the two subsets of type $A_5$ containing $N_2$,  $D_4$, the subset of type $4A_1$ which is stable under the automorphism of $\tilde{\Delta}$, $\{N_0,\,N_2,\,N_3\}$ and $\{N_2,\,N_5,\,N_7\}$;
  \item[$E_8$ :] $\tilde{\Delta}\setminus\{N_1,N_3\}$, $\tilde{\Delta}\setminus\{N_1,\,N_3,\,N_6\}$, $\tilde{\Delta}\setminus\{N_4,N_6,\,N_8\}$, 
  $\{N_2,N_5,N_7,N_0\}$  or of type $D_7$, $E_7$, $D_6+A_1$, $2A_3+A_1$, $3A_2+A_1$, $D_5+2A_1$, $D_4+A_3$, $D_6$, $E_6$, $D_4+2A_1$, $3A_2$, $D_4$; 
 \item[$F_4$ :] of type $A_3$, $A_1+B_2$, $2A_1+\tilde{A}_1$, $B_3$, $C_3$,
  $2A_1$, $B_2$, $\tilde{A}_2$;
  \item[$G_2$ :] of type $\tilde{A}_1$.
  \end{itemize}
  If $S_F$ is in this list, then $X$ is also unibranch.
\end{proposition}
\pf If $S_F=\emptyset$, then $X=T/W$, whereas if $\#S_F=\#\Delta$, then $X$ is a point and there is nothing to prove. For the remaining cases, we know from Lemma \ref{lem:effe} that $X$ is normal in codimension $1$ if and only if $\#\Sigma(\HH,L)=1$, i.e., if and only if $[S_F]=\{S_F\}$. By Corollary \ref{cor:algo} from which we adopt notation, this happens if and only if for every $j\in\hat{S}_F$ we have $w_{0j}S_F=S_F$. By inspection we obtain the given list. 
\epf

We are ready to produce the full list of normal and smooth strata.

\begin{theorem}\label{thm:lista-normali}Let $X=X(\HH,L)$ be a stratum in $\Ec/W_{aff}$ and let $F$ be a maximal face in $\Sigma(\HH,L)$.
If $\Phi$ is classical, then $X$ is normal if and only if it is normal in codimension 1. If $\Phi$ is exceptional, then
$X$ is normal, or equivalently smooth, if and only if $S_F=\emptyset$, or $\#S_F=\#{\Delta}$, or it is as follows:
\begin{itemize}
 \item[$E_6$ :]  of type $A_5$ (there are three such subsets), $4A_1$, $2A_2$ (there are three such subsets);
  \item[$E_7$ :] of type $E_6$, $D_6$ (there are two such subsets), $D_5+ A_1$ (there are two such subsets), $D_4+2A_1$, $2A_3$, $3A_2$, $A_3+3A_1$ (there are two such subsets),  $5A_1$, the two subsets of type $A_5$ containing $N_2$,  the subset of type $4A_1$ which is stable under the automorphism of $\tilde{\Delta}$, $\{N_0,\,N_2,\,N_3\}$ and $\{N_2,\,N_5,\,N_7\}$;
  \item[$E_8$ :] $\tilde{\Delta}\setminus\{N_1,N_3\}$, $\tilde{\Delta}\setminus\{N_1,N_3,\,N_6\}$, $\tilde{\Delta}\setminus\{N_4,N_6,\,N_8\}$, 
  $\{N_2,N_5,N_7,N_0\}$  or of type $D_7$, $E_7$, $D_6+A_1$, $2A_3+A_1$, $3A_2+A_1$, $D_5+2A_1$, $D_4+A_3$, $3A_2$; 
 \item[$F_4$ :] of type $A_3$, $A_1+B_2$, $2A_1+\tilde{A}_1$, $B_3$, $C_3$,
  $2A_1$, $\tilde{A}_2$;
  \item[$G_2$ :] of type $\tilde{A}_1$.
  \end{itemize}
  \end{theorem}		
\pf We only need to consider the strata listed in Proposition \ref{prop:codim1}. If $S_F=\emptyset$, then $X=T/W$, whereas if $\#S_F=\#\Delta$, then $X$ is a point and there is nothing to prove. If $\#S_F=\#\Delta-1$ then $\dim L=1$ so $X$ is normal because it is normal in codimension $1$. In the remaining cases, we observe that $X$ is unibranch by Lemma \ref{lem:effe}. Lemma \ref{lem:unions} and Proposition \ref{prop:constant} imply that normality should be verified at points in $\overline{F}$. Here Corollary \ref{cor:X'} applies and by Theorem \ref{thm:isolated1} it is enough to check normality of $X(\HH_l,L)$ for $l$ ranging in the set of vertices of $\overline{F}$, i.e., in the set of  nodes in $\hat{S}_F$. Let $N_l$ be the node corresponding to vertex $l$. The stabiliser $\WHl$ is generated by the reflections with respect to all hyperplanes containing $l$, i.e., by the reflections corresponding to all nodes but $N_l$. Its Coxeter graph is thus obtained from the Coxeter graph of $\WH$ by removing $N_l$. The Coexeter class corresponding to $\Sigma(\HH_l,L)$ contains the unique subset $S_{F_l}=S_F$ by locality of the algorithm in Corollary \ref{cor:algo}. Also, $S_{F_l}$ is obtained by removing $N_l$ from the graph in $S_F$. 
The parametrization in terms of subsets of the Coxeter graph coincides with the one used in \cite{broer, DR, richardson}. In other words, $X(\HH_l,L)$ is normal if and only if the subset $S_F$ occurs in $\cite[Tables\  I, II]{DR}$ for the Coxeter group whose generating system is obtained by removing the node $N_l$ from $S$. The required list is obtained by checking this property for all nodes in $\hat{S}_F$.
\epf

\section{The general case}\label{sec:non-sc}

In this Section we give some criteria to analyse the case of general $K$ and provide some key examples of strata for $K\neq 1$. These results will be applied in Section \ref{sec:simple}
to  obtain the complete list of normal strata for the Jordan stratification in  $G/\!/G$ for any simple $G$.

\subsection{Relative criteria}
In this Section we provide criteria for unibranchedness and for normality of $X=X(\HH, K,L)$ in terms of $X'=X(\HH, K', L)$  for $K'\lhd K $ two admissible subgroups.
%of  automorphisms of $\Abar$ as in Section \ref{subsec:basic}. 
These criteria will be applied in Section \ref{sec:non-sc} to the special case of $\WH$ an irreducible affine Weyl group and $K,\,K'$ subgroups of $P^\vee/Q^\vee$, i.e., for Jordan strata in isogenous groups.

If $K'\lhd K$ are as above, then $I:=K/K'$ \label{I} acts on $E_{\mathbb{C}}/\WW_{\HH,K'}$. Let $\pi : E_{\mathbb{C}}/\WW_{\HH,K'} \rightarrow E_{\mathbb{C}}/\WK$ be the canonical projection. Clearly, $\pi^{-1}(X)=IX'$ and $X\simeq IX'/I$.

Let $I^{X'}$ be the stabilizer of $X'$ in $I$ and, for  $x\in X'$,  let $I^{X',x}:=I^{X'}\cap I^x$.  By abuse of notation we will indicate points in a quotient by a representative.\label{inertia}
\begin{lemma}
Let $X=X(\HH,K,L)$ and $X'=X(\HH,K',L)$ with $K'\lhd K$ and $I=K/K'$.
\begin{enumerate}[label=(\roman*)]
%\item We have $\pi^{-1}(X)=I X'$ and $\pi$ induces an isomorphism $\pi^{-1}(X)/I\simeq X$.
\item\label{boh}
The stratum  $X$ is unibranch at the point $x$ if and only if $I^x X'/ I^x$ is unibranch at $x$ and 
$\left(I^x X',x\right)\sim_{se}\left(\pi^{-1}(X),\,x\right)$.
\item
If \ref{boh} holds, then $X$ is normal at $x$ if and only if $X'/ I^{X',x}$ is normal at $x$ and the canonical map $X' /I^{X',x}\rightarrow X$ is a local analytic isomorphism around $x$.
\end{enumerate}
\end{lemma}
\pf
\begin{enumerate}[label=(\roman*)]
%\item Follows from the definition.
\item Let $Y_x$ be the union of the irreducible components of $I X'$ containing $x$. There is a small enough $I$-stable analytic neighbourhood $U_x$ of $x$ such that $Y_x\cap U_x= I X'\cap U_x$. Combining with Lemma \ref{lem:invariants}, we obtain $\left(X,x\right)=\left(I X'/ I,\,x\right)\sim_{se}\left(Y_x/ I^x,\,x\right)$. Since $I^x X'/ I^x$ is an irreducible component of $Y_x / I^x$, the latter is unibranch if and only if  $I^x X' / I^x = Y_x / I^x$ is unibranch at $x$, whence also $I^xX'\cap U_x=IX'\cap U_x$. 
\item By the proof of \ref{boh} normality of $X$ at $x$ is equivalent to normality of $I^x X' / I^x$ at $x$. The map $f: X'/ I^{X',x} \rightarrow I^x X'/ I^x$ induced by the inclusion of $X'$ in $I^x X'$ is a finite birational morphism, hence the source and the target varieties have the same normalisation. So, the target is normal at $x$ if and only if the source is normal at $x$ and  $f$ is an isomorphism at $x$.
\end{enumerate}
\epf

We will mainly use the following special case. 
\begin{corollary}\label{cor:soprasotto} 
If $I^{X',x} = 1$, then $X(\HH,K,L)$ is normal at $x$ if and only if it is unibranch at $x$, $X(\HH,K',L)$ is normal at $x$, and $(X(\HH,K,L),x)\sim_{se}(X(\HH,K',L),x)$. In the special case of $I^x = 1$, the first two conditions suffice.
\end{corollary}

\subsection{Some examples}
In this subsection $\HH$ is effective and finite, so $\WH$ fixes a point which we set as an origin $O$. Hence the Euclidean space $E$ is identified with its vector space of translations $V$. 

\begin{example}\label{ex:dnbn}
{\rm Let $\WH=W$ be the Weyl group of type $D_n$ and $\Abar$ be the fundamental chamber. It is a simplicial cone with vertex $0$ generated by the half-lines ${\mathbb R_{\geq0}}\omega_i^\vee$ for $i\leq n$. Let $K$ be generated by the orthogonal transformation $k$ fixing $\omega_i^\vee$ for $i\leq n-2$ and interchanging $\omega_{n-1}^\vee$ and $\omega_n^\vee$. Then, $\WK$ is again a finite reflection group $\WW_{\HH_B}$ for  the central hyperplane arragement  $\HH_B$ obtained by adding to $\HH$ the hyperplane containing every $\omega_i^\vee$ for $i\neq n-1$ and $\omega_{n-1}'=\omega_{n-1}^\vee+\omega_n^\vee$. It is the Weyl group of type $B_n$ and its fundamental chamber $\Abar_B$ is the simplicial cone generated by the half-lines ${\mathbb R_{\geq0}}\omega_i^\vee$ for $i\neq n-1$ and $R_{\geq0}\omega_{n-1}'$.

A flat $L$ for $\HH$ lying over $\Abar$ lies also over $\Abar_B$ provided that if $F=\interior\cap\Abar$ has ${\mathbb R_{\geq0}}\omega_{n-1}^\vee$ as generating line, then it also has ${\mathbb R_{\geq0}}\omega_{n}^\vee$ as generating line. In terms of nodes, it means that if $N_{n-1}\in\hat{S}_F$, then $N_{n}\in\hat{S}_F$. So if $L$ lies over $\Abar$, at least one flat among $L$ and $kL$ lies over $\Abar_B$ and $\Abar$. 

If $L$ lies over $\Abar$ and $\Abar_B$, then we have equality of strata $X(\HH,K,L)=X(\HH_B,L)$. 

In this situation, let $F=\interiorc\cap\Abar$, $F'=\interiorc\cap\Abar_B$. If $F$ is generated by lines with indices $\neq n-1$, then $F'=F$. If $F$ has ${\mathbb R_{\geq0}}\omega_{n-1}^\vee$ and ${\mathbb R_{\geq0}}\omega_{n}^\vee$ as generating lines, then $F'$ has ${\mathbb R_{\geq0}}\omega_{n}^\vee$ and ${\mathbb R_{\geq0}}\omega_{n-1}'$ as generating lines. In both cases, $S_F$ and $S_{F'}$ contain the nodes corresponding to the same index set. 

By \cite[Proposition 8.2.1]{richardson} the normal strata in $\Ec/\WK$ are those for which either $S_F=\emptyset$, or $S_F=\Delta$, or it is of type $D_{m_0}+ dA_h$ with $n=m_0+n_0+d(h+1)$, $h\geq0$, $m_0\geq0$. } 
\end{example}
\begin{example}\label{ex:giro}{\rm
Let $\WH\simeq \prod_{j=1}^t\WW_{\HH_{(j)}}$ with $W_{\HH_{(j)}}\simeq\WW_{\HH_{(i)}}$ for any $j,\,i$ act on $E=\prod_{j=1}^tE_{(j)}$ and let $K=\langle k\rangle \simeq {\mathbb Z}/t{\mathbb Z}$ act permuting components of $E$ cyclically.
Assume that $L=L_{(1)}\times \prod_{j=2}^t\{p_j\}$ has trivial component on $E_{(j)}$ for $j\neq 1$.
%and that $X(\HH_{(1)},1,L_{(1)})$
%= \WW_{\HH_{(1)}} L_{(1)\mathbb C}/\WW_{\HH_{(1)}}$ 
%is normal. 
Then,
%$\Gamma_{\HH,K,L}=\Gamma_{\HH,L}=\Gamma_{\HH_{(1)}, L_{(1)}}\times\left(\prod_{j=2}^t\WW_{\HH_{(j)}}\right)$ and 
\begin{align*}
&{\mathbb C}[\Ec]^{\WK}\simeq\left(\otimes_{j=1}^t{\mathbb C}[E_{(j)}]^{\WW_{\HH_{(j)}}}\right)^K\simeq {\mathbb C}[E_{(1)}]^{\WW_{\HH_{(1)}}},\\
&\Gamma_{\HH,K,L}=\Gamma_{\HH,L}=\Gamma_{\HH_{(1)}, L_{(1)}}\times\left(\prod_{j=2}^t\WW_{\HH_{(j)}}\right),\\
&{\mathbb C}[\Lc]^{\Gamma_{\HH,K,L}}\simeq {\mathbb C}[L_{(1),{\mathbb C}}]^{\Gamma_{\HH_{(1)},L_{(1)}}}.
\end{align*}
 and by Proposition \ref{prop:richardson} $X(\HH,K,L)$ is normal if and only if 
$X(\HH_{(1)},L_{(1)})$ is so.}
%
%Then  $X(\HH,K,L)\simeq \bigcup_{j=1}^t(\prod_{l=j+1}^t\{p_l\}\times X(\HH_{(1)},1,L_{(1)})\times\prod_{l=2}^{j}\{p_l\})/K$ is isomorphic to the image of $X(\HH_{(1)},1,K)$ in $(X(\HH_{(1)},1,L_{(1)}))
%^t/K$. Therefore, $X(\HH,K,L)\simeq X(\HH_{(1)},1,L_{(1)})\simeq X({\HH},1,L)
%$ is normal. }
\end{example}

\subsection{Simple groups}\label{sec:simple}

In this Section we deal with Jordan strata in simple groups, i.e.,  $\WH=W_{aff}=W\ltimes Q^\vee$ is irreducible and $K\leq P^\vee/Q^\vee$. Here $\Abar$ is the fundamental alcove.

Observe that if a stratum $X(\HH,K,L)$ is normal then there is $F\in\Sigma(\HH,L)$
such that  $\Sigma(\HH,K,L)=\{kF~:~k\in K\}$ by Lemma \ref{lem:condition}. Also, Proposition \ref{prop:nec_codim1} and Lemma \ref{lem:idd} imply that \eqref{eq:condition} holds for any $L'\subset L$ lying over $\Abar$. 
%by Proposition \ref{prop:equivalent} and Lemma \ref{lem:idd}. 
These combinatorial conditions can be verified by looking at the action of $K$ on the set of vertices of $\Abar$ or, equivalently, on the corresponding set of nodes in $S$. A vertex $x_j$ lies in the closure of a face in $\Sigma(\HH,K,L)$ if and only if $N_j\in\hat{S}_F$.
In particular,  if a vertex $x_j\in\overline{F}$ has trivial stabiliser in $K$, taking $L'=\{x_j\}$ in \eqref{eq:condition} together with $\Sigma(\HH,K,L)=\{kF~:~k\in K\}$ gives the necessary condition for normality of $X(\HH,K,L)$:
\begin{equation}\label{eq:nodes}
%(K^{x_j}=1)\implies \{k\hat{J}_F~:~k\in K,\, N_j\in k\hat{S}_F\}=\{\hat{S}_F\}.
\textrm{If } K^{x_j}=1\textrm{ and }N_j\in \hat{S}_F\cap k\hat{S}_F \textrm{ for some }k\in K\implies k\hat{S}_F=\hat{S}_F.
\end{equation}

When $K$ is small we also have the following necessary condition.
\begin{lemma}\label{lem:allbutone}Assume $\#K=2$ and let $X=X(\HH,K,L)$ be the stratum corresponding to $S_F\subset S$. If $X$ is normal, then either $X(\HH,L)$ is normal in codimension $1$ or else $\hat{S}_F$ has exactly one vertex that is not fixed by $K$ and $X(\HH,L)$ is unibranch. 
\end{lemma}
\pf If $X$ is normal, then either  $\Sigma(\HH,L)=\{F\}$, and Lemma \ref{lem:effe} applies, or else $\Sigma(\HH,L)=\Sigma(\HH,K,L)=\{F,\,kF\}$, where $k$ is the non-trivial element in $K$. In this case, $X(\HH,L)$ is unibranch by Lemma \ref{lem:due-uni}. Also, Corollary \ref{cor:connected} shows that the faces $F$ and $kF$ must be separated by a wall, so $\hat{S}_F\cap k\hat{S}_F$ contains all nodes of $\hat{S}_F$ but $1$. By \eqref{eq:nodes}, all such nodes are fixed by $K$.   
\epf

The following  special case can be treated directly.

\begin{lemma}\label{lem:ferma}Let $X=X(\HH,K,L)$ be a stratum and let $F\in \Sigma(\HH,K,L)$ be such that $kF=F$ for any $k\in K$. Then
\begin{enumerate}[label=(\roman*)]
\item If $X$ is normal in codimension $1$, then $X(\HH,L)$ is normal in codimension $1$ and $X$ is unibranch.
\item If $X(\HH,L)$ is normal in codimension $1$, then $X$ is normal in codimension $1$ and unibranch.  
\item If $X(\HH,L)$ is normal then $X$ is normal.
\item If $W$ is classical, then $X(\HH,L)$ is normal if and only if $X$ is normal.
\end{enumerate}
\end{lemma}
\pf  (i). By Lemma \ref{lem:condition} we necessarlly have $\Sigma(\HH,K,L)=\{F\}=\Sigma(\HH,F)$, so $X(\HH,L)$ in normal in codimension $1$ by Lemma \ref{lem:effe}. Also, \eqref{eq:condition} holds for every flat, so Lemmata %\ref{prop:equivalent} and 
%Lemma 
\ref{lem:unions} and \ref{lem:id} imply  that $X$ is unibranch.

(ii) If $X(\HH,L)$ is normal in codimension $1$, then $\Sigma(\HH,K,L)=K\Sigma(\HH,L)=\{F\}$. Hence, \eqref{eq:condition} holds for every flat $L'$ lying over $\Abar$, so Proposition \ref{prop:equivalent} and 
Lemma \ref{lem:id} imply  that $X$ is normal in codimension $1$ and unibranch.

(iii) If $X(\HH,L)$ is normal,  then $X=K X(\HH,L)/K=X(\HH,L)/K$ is normal.

(iv) One direction is (iii). The other follows from (i) because for classical $W$ normality of and normality in codimension $1$ coincide for $X(\HH,L)$, by Theorem \ref{thm:lista-normali}.
\epf

We will deal with each irreducible root system and non-trivial choice of $K\leq P^\vee/Q^\vee$ separately. We recall that numbering of simple reflections and nodes in the Coxeter graph of $\WH$ are as in \cite{bourbaki}. 
We will constantly make use of the identification in Proposition \ref{prop:stratifications}. Recall that Remark \ref{rem:real}, Theorem \ref{thm:isolated1} and Lemma \ref{lem:unions} imply that it is enough to verify normality at all vertices of faces in $\Sigma(\HH,K,L)$.

\subsubsection{\bf Type $A_n$}

In this case $P^\vee/Q^\vee$ permutes cyclically the nodes in $S$, whence $K^l=1$ for any vertex $l$ of $\Abar$ and any choice of $K$. We set $K=\langle k\rangle$.

\begin{proposition}
Let $G$ be a group of type $A_n$. A stratum $X=X(\HH,K,L)$ in the Jordan stratification of $G/\!/G$ is normal if and only if the stratum $X(\HH,L)$ in $\mathrm{SL}_{n+1}/\!/\mathrm{SL}_{n+1}$ is so.
\end{proposition}
\pf If $X$ is normal, then $X(\HH,L)$ is normal at each vertex $x_j$ of any face in $\Sigma(\HH,K)$ by Corollary \ref{cor:soprasotto}, so $X(\HH,L)$ is normal. 
%Let $F\in\Sigma(\HH,L)$. If $kF=F$ then we apply Lemma \ref{lem:ferma} (iv).
%Assume $kF\neq F$ and that $X$ is normal. Then 
% %By Theorem \ref{thm:lista-normali} and Lemma \ref{lem:effe} it is enough to show that 
% $|\Sigma(\HH,L)|=1$. Were this not the case, there would be two maximal faces $F_1,\,F_2$ in $\Sigma(\HH,L)$ that are separated by a wall $H$ of $\Abar$ by Corollary \ref{cor:connected}. But then, there would be at least one vertex $x_j$ of $\Abar$ lying in $\overline{F_1}\cap\overline{F_2}$. This would contradict unibranchedness of  
% $X$ at $x_j$ in virtue of Proposition \ref{prop:equivalent}, condition \eqref{eq:condition} and the equality $K^{x_j}=1$. 
 
If $X(\HH,L)$ is normal, then $\Sigma(\HH,L)=\{F\}$ % and $\WW_{\HH^L}\leq \iota_L(\Gamma_{\HH,L})\leq  \iota_L(\Gamma_{\HH,K,L})$
 by Proposition \ref{prop:equivalent}. Using the list of normal strata in Theorem \ref{thm:lista-normali} we verify that in type $A_n$ there holds either $\hat{S}_{kF}=k\hat{S}_F=\hat{S}_F$ or else
$\hat{S}_{kF}\cap \hat{S}_F=\emptyset$, i.e., either $kF=F$ or else $k\overline{F}\cap \overline{F}=\emptyset$. In the first case Lemma \ref{lem:ferma} applies. In the second case, 
$\Sigma(\HH,K,L)=\{kF~:~k\in K\}$ and $X$ is unibranch at all vertices of faces in $\Sigma(\HH,K,L)$ by Lemma \ref{lem:set}.  Corollary \ref{cor:soprasotto} gives normality of $X$.\epf

\subsubsection{\bf Type $B_n$ for $n\geq3$}

Here $K=P^\vee/Q^\vee\simeq\mathbb{Z}/2\mathbb{Z}$, corresponding to $G=\mathrm{SO}_{2n+1}$. The non-trivial element $k\in K$ acts on the vertices of $\Abar$ interchanging the vertices $x_0=0$ and $x_1=\omega_1^\vee$, fixing their middle point $x_1'=\frac{1}{2}\omega_1^\vee$ and $x_j$ for $j=2,\ldots,\,n$. In this case, $\WK=(Q^\vee\rtimes W)\rtimes K$ is again a reflection group. It is the group $\WW_{\HH_C}$ for the affine hyperplane arrangement $\HH_C$ obtained by adding to $\HH$ all $\WH$-translates of the affine hyperplane $H'$ passing through $x_1'$  and $x_j$ for $j=2,\ldots,\,n$. In other words, it is 
%Thus, $\WK$ is
%isomorphic to 
the affine Weyl group of type $C_n$.

The fundamental alcove $\Abar_C$ for $\WK$ has vertices $x_0$, $x_1'$, and $x_j$ for $j\geq2$ and its walls are $H'$ and the walls of $\Abar$ except from the hyperplane $H$ containing the vertices $x_j$ for $j=1,\,\ldots,\,n$, i.e., the wall corresponding to the node labeled by $0$. We denote by $S'$ the set of reflections with respect to these walls and we identify it with the set of nodes in the Coxeter graph of $\WW_{\HH_C}$. 

\begin{lemma}\label{lem:translation}
Let $L$ be a flat for $\HH$ and let $X=X(\HH,K,L)$. Then
\begin{enumerate}[label=(\roman*)]
 \item There is always a $\WK$-conjugate of $L$ lying over $\Abar$ and $\Abar_C$.
 \item Assume $L$  lies over $\Abar$ and $\Abar_C$ and let $F=\interior\cap\Abar$, $F'=\interior\cap\Abar_C$. Then $X$ is isomorphic to the stratum in $\Ec/\WW_{\HH_C}$ indexed by the subset $S_{F'}$ of $S'$ consisting of the nodes with same indices as $S_F$.
\end{enumerate}
\end{lemma}
\pf (i) We may always assume that $L$ lies over $\Abar$. Since $L$ is also a flat for $\HH_C$, there is $w\in \WW_{\HH_C}=\WK$ such that  $w L$ lies over $\Abar_C$. Since $\Abar_C\subset\Abar$, by Remark \ref{rem:fundamental} there is $k\in K$ such that $wL=kL$, hence $wL$ lies over $\Abar$ and $\Abar_C$. 

(ii) Observe that a flat $L$ lying over $\Abar$ lies also over $\Abar_C$ unless $x_0\not\in\overline{F}$ and $x_1\in\overline{F}$. Therefore, if $L$ lies over $\Abar$ and $\Abar_C$ and $x_1\not\in\overline{F}$, then $F=F'$. If $x_1\in\overline{F}$, then $x_0\in\overline{F}$ so $x_1'\in\overline{F}$ and the vertices of  $F'$ are obtained from the vertices of $F$ by replacing $x_1$ by $x_1'$. Hence the indices involved in $S_F$ and $S_{F'}$ coincide. 
\epf

\begin{proposition}
Let $X=X(\HH,K,L)$ be a stratum in the Jordan stratification of $\mathrm{SO}_{2n+1}/\!/\mathrm{SO}_{2n+1}$ corresponding to the subset $S_F$  of the Coxeter graph of $\WH$. Then $X$ is normal if and only if $S_F=\emptyset$, or $\#S_F=\#\Delta$ or it is of type $D_{m_0}+ dA_h + B_{n_0}$ with $n=m_0+n_0+d(h+1)$ and either $m_0\geq 1$, $n_0\geq 0$, $h\geq 0$, or else $n_0=0$, $m_0\geq0$ and $h\geq0$. In particular, a stratum is normal if and only if it is smooth.
\end{proposition}
\pf This is obtained applying Lemma \ref{lem:translation} to Theorem \ref{thm:lista-normali}. Last statement holds because of the identification with strata in $\Ec/\WW_{\HH_C}$.  
\epf

Comparing with Theorem \ref{thm:lista-normali} we see that there are strata $X(\HH,L)$ in $\mathrm{Spin}_{2n+1}/\!/\mathrm{Spin}_{2n+1}$ that are not normal even if the corresponding strata $X(\HH,K,L)$ in  $\mathrm{SO}_{2n+1}/\!/\mathrm{SO}_{2n+1}$ are so.

\subsubsection{\bf Type $C_n$ for $n\geq2$}

Here $K=P^\vee/Q^\vee\simeq\mathbb{Z}/2\mathbb{Z}$, corresponding to $G=\mathrm{PSp}_{2n}$. The non-trivial element $k\in K$ on $\Abar$ interchanges the vertices $x_j$ and $x_{n-j}$ in $\Abar$ for $j=0,\,\ldots,\,n$. If $n=2m$ is even, then $kx_m=x_m$, whereas if $n$ is odd, $K^x=1$ for every vertex of $\Abar$.

\begin{proposition}
Let $X=X(\HH,K,L)$ be the stratum in the Jordan stratification of $\mathrm{PSp}_{2n}/\!/\mathrm{PSp}_{2n}$ corresponding to the Coxeter class $[S_F]$.
Then $X$ is normal if and only if $X(\HH,L)$ is normal and \eqref{eq:nodes} holds.
\end{proposition}
\pf If $k\hat{S}_F=\hat{S}_F$, then if either $X$ or $X(\HH,L)$ is normal \eqref{eq:nodes} automatically holds and the statement follows from  Lemma \ref{lem:ferma}. We assume for the rest of the proof that $k\hat{S}_F\neq\hat{S}_F$. 
If $X$ is normal, then \eqref{eq:nodes} holds. %by Lemma \ref{lem:idd}. 
Since there is at most one node in $\hat{S}_F$ that is fixed by $K$, Lemma \ref{lem:allbutone} shows that the only possibility for $X(\HH,L)$ not being normal could be for $\#\hat{S}_F\leq2$ with possible equality only for $n=2m$ and $x_m\in\overline{F}$. However, in this situation $X(\HH,L)$ is always normal by Theorem \ref{thm:lista-normali}.  

Conversely, assume that $X(\HH,L)$ is normal and \eqref{eq:nodes} holds. Then 
%$\WW_{\HH^L}\leq \iota_L(\Gamma_{\HH,K,L})$ so  
Lemma \ref{lem:id} gives unibranchedness at all vertices of $F$ and $kF$. By Corollary \ref{cor:soprasotto} the stratum $X$ is normal at all vertices of $F$ and $kF$ with trivial stabiliser. This concludes the discussion for $n$ odd. Assume now $n=2m$ and $x_m\in\overline{F}$. By Corollary \ref{cor:X'} it is enough to prove normality of the finite counterpart $X(\HH_{x_m},K,L)$ at $x_m$.

The finite counterpart  $X(\HH_{x_m},L)$ of $X(\HH,L)$ at $x_m$ is a Jordan stratum for the semisimple Lie algebra of type $C_m\times C_m$.  By inspection we see that, up to $K$-action, the subsets $S_F$ from Proposition \ref{prop:codim1} satisfying condition \eqref{eq:nodes} and not containing $N_m$ necessarily contain all nodes with indices $\geq m$, so  $X(\HH_{x_m},L)\simeq X(\HH_{x_m(1)},L_{(1)})\times\{p\}$
where $X(\HH_{x_m(1)},L_{(1)})$ is a normal stratum for the finite Coxeter group of type $C_m$ and $K$ interchanges the two components of $\hh$. We are in the situation of Example \ref{ex:giro}, so  $X(\HH_{x_m},K,L)$ is normal.  
\epf

\subsubsection{\bf Type $D_n$ for $n\geq 4$}

If $n$ is odd, $P^\vee/Q^\vee\simeq \mathbb{Z}/4\mathbb{Z}\simeq\langle \sigma\rangle$ where the action of $\sigma$ on the vertices of $\Abar$ is given by $x_0\mapsto x_n\mapsto x_1\mapsto x_{n-1}\mapsto x_0$ and $x_{j}\mapsto x_{n-j}$ for $2\leq j\leq n-2$.

\noindent If $n$ is even, $P^\vee/Q^\vee\simeq \mathbb{Z}/2\mathbb{Z}\times \mathbb{Z}/2\mathbb{Z}\simeq\langle \tau_1\rangle\times\langle \tau_2 \rangle$, where the action of 
$\tau_1$ is given by 
$x_0\mapsto x_1\mapsto x_0$, $x_n\mapsto x_{n-1}\mapsto x_n$ and $x_{j}\mapsto x_{j}$ for $2\leq j\leq n-2$ and the action of 
$\tau_2$ is given by 
%$x_0\mapsto x_n\mapsto x_0$, $x_1\mapsto x_{n-1}\mapsto x_1$ and 
$x_{j}\mapsto x_{n-j}$ for $0\leq j\leq n$.

\bigskip 

Let us consider the case in which $K=\langle \xi\rangle$ is the group of order $2$ with $\xi=\tau_1$ when $n$ is even and $\xi=\sigma^2$ when $n$ is odd. The stratification in $\Ec/\WK$ is the Jordan stratification in $\mathrm{SO}_{2n}/\!/\mathrm{SO}_{2n}$.
 
\begin{proposition}\label{prop:so}
Let $X=X(\HH,K,L)$ be the stratum in the Jordan stratification of $\mathrm{SO}_{2n}/\!/\mathrm{SO}_{2n}$ corresponding to the Coxeter class $[S_F]$. 
Then $X$ is normal if and only if one of the following two conditions hold:
\begin{enumerate}[label=(\roman*)]
\item $X(\HH,L)$ is normal;
\item $S_F$ is of type $D_{m_0}+dA_h$ with  $n=m_0+d(h+1)$, $m_0\geq2$ and $h$ is even.
\end{enumerate}
\end{proposition}
\pf If $\xi{S}_F={S}_F$, then (ii) cannot occur and  Lemma \ref{lem:ferma} gives equivalence of normality of $X$ with (i).  We suppose from the rest of the proof that $\xi S_F\neq S_F$. If $X$ is normal, then Lemma \ref{lem:allbutone} implies that either $X(\HH,L)$ is normal and $\Sigma(\HH,K,L)=\{F,\,\xi F\}$, or else
$\#\hat{S}_F\cap \{N_0,N_1,N_{n-1},N_n\}=1$  and $\Sigma(\HH,L)=\Sigma(\HH,K,L)=\{F,\,\xi F\}$.  In this case, the algorithm in Corollary \ref{cor:algo} shows that $S_F$ is necessarily as in (ii). 

Conversely, assume that either (i) or (ii) hold. Then, $\Sigma(\HH,K,L)=\{F, \xi F\}$ and $X(\HH,K,L)$ is unibranch at all vertices of $F$ and $\xi F$ by Lemma \ref{lem:id} and $X(\HH,L)$ is unibranch by Lemma \ref{lem:due-uni}.
Thus, $X(\HH,L)$ is normal at the vertex of $\overline{F}$ with trivial stabiliser in case (ii) by Corollary \ref{cor:X'} and  \cite[8.3.1]{richardson}, whence Corollary \ref{cor:soprasotto} gives normality at all vertices of $\overline{F}$ and $\overline{\xi F}$ with trivial stabiliser. 
By Corollary \ref{cor:X'} it remains to prove that  $X(\HH_{x_j},K,L)$ is normal for every $x_j\in\overline{F}$ with $2\leq j\leq n-2$.

These are strata for  a finite Coxeter group of type $D_j\times D_{n-j}$. If $S_F$ has a component of type $D_{m_0}$ for $m_0\geq 2$, then  $X(\HH_{x_j},K,L)\simeq X(\HH_{x_j(1)},K,L_{(1)})\times X(\HH_{x_j(2)},K,L_{(2)})$ where $K$ acts trivially on $L_{(1)}$.  Thus
$X(\HH_{x_j(1)},K,L_{(1)})=X(\HH_{x_j(1)},L_{(1)})$  is normal by \cite[8.3.1]{richardson}, whereas  $X(\HH_{x_j(2)},K,L_{(2)})$  is isomorphic to a normal stratum for a Coxeter group of type $B_{n-j}$ in virtue of Example \ref{ex:dnbn}. If, instead, $S_F$ is of type $dA_h$, we consider the extension $\tilde{K}=\langle K,\xi'\rangle$ of $K$ where $\xi'$ is the involution of $\Abar$ interchanging $x_0$ and $x_1$ and fixing all other veritces. Then, 
$X(\HH_{x_j},\tilde{K},L)=\left(\langle\xi'\rangle X(\HH_{x_j},K,L)\right)/ \langle\xi'\rangle$ is normal  in virtue of Example \ref{ex:dnbn}, so $X(\HH_{x_j},K,L)$ is normal by Corollary \ref{cor:soprasotto} applied to $K\lhd \tilde{K}$.
\epf

Let us now consider $K=\langle \tau_2\rangle$ for $n=2m$ even, so $K$ fixes only the vertex $x_m$. The stratification in $\Ec/\WK$ is the Jordan stratification in $\mathrm{HSpin}_{2n}/\!/\mathrm{HSpin}_{2n}$.

\begin{proposition}
Let $X=X(\HH,K,L)$ be the stratum in the Jordan stratification of $\mathrm{HSpin}_{2n}/\!/\mathrm{HSpin}_{2n}$ corresponding to the Coxeter class $[S_F]$. 
Then $X$ is normal if and only if the stratum  $X(\HH,L)$ in $\mathrm{Spin}_{2n}/\!/\mathrm{Spin}_{2n}$ is normal and \eqref{eq:nodes} holds.
\end{proposition}
\pf If $\tau_2{S}_F={S}_F$,  Lemma \ref{lem:ferma} (iv) shows that $X(\HH,L)$ is normal if and only if $X$ is so and if this is the case, \eqref{eq:nodes} automatically holds.
 %and  $X(\HH,L)$ is normal if and only if $X$ is so. 
We assume for the rest of the proof that $\tau_2{S}_F\neq{S}_F$ and that $X$ is not a point, i.e., $\#\hat{S}_F\geq2$. 
 
If $X$ is normal, then \eqref{eq:nodes} holds and Lemma \ref{lem:allbutone} shows that either $X(\HH,L)$ is normal, or else $\#\hat{S}_F=2$ and $S_F$ is of type $D_{m}+A_{m-1}$ with $m$ odd. However, the latter case is ruled out because Corollary \ref{cor:algo} would yield  $\#\Sigma(\HH,K,L)\geq 4$ contradicting Lemma \ref{lem:condition}. 

Assume now $X(\HH,L)$ is normal and \eqref{eq:nodes} holds.  Then $\Sigma(\HH,K,L)=\{F,kF\}$ and 
%\eqref{eq:nodes} together with
 Lemma \ref{lem:id} guarantees unibranchedness at all vertices of $\overline{F}$ and $k\overline{F}$. Normality of $X$ at all vertices but $x_m$ follows from Corollary \ref{cor:soprasotto}. We focus on $x_m$. Observe that if $N_m\in\hat{S}_F$ and $S_F\neq kS_F$ with $S_F$ of type $dA_h$, then it does not satisfy  \eqref{eq:nodes}. If it is of type $D_{m_0}+D_{n_0}+dA_h$, for $m_0>n_0\geq0$ as in Proposition \ref{prop:codim1}, then \eqref{eq:nodes} 
forces $m_0=m> n_0\geq0$. Hence, $X(\HH_{x_m},L)$ has two factors one of which is trivial whereas the other is normal, and $K$ interchanges the two factors. Thus, $X$ is normal at $x_m$ as in Example \ref{ex:giro}.
\epf

Finally, we consider the group $K=P^\vee/Q^\vee$. The stratification in $\Ec/\WK$ is the Jordan stratification in $\mathrm{PSO}_{2n}/\!/\mathrm{PSO}_{2n}$. We write $\eta$ for either $\sigma$ or $\tau_2$ and $\xi$ for either $\sigma^2$ or $\tau_1$. Recall that $\tau_2$ fixes only $x_m$ and that $K^{x_j}=1$ holds only for $j\in\{0,1,n-1,n\}$.

\begin{proposition}
Let $X=X(\HH,K,L)$ be the stratum in the Jordan stratification of $\mathrm{PSO}_{2n}/\!/\mathrm{PSO}_{2n}$ corresponding to the Coxeter class $[S_F]$. 
Then $X$ is normal if and only if the following two conditions hold:
\begin{enumerate}[label=(\roman*)]
\item $\{k \hat{S}_F\,:\,k\in K^{x_j}\}=\{k\hat{S}_{F_1}\,:\,S_{F_1}\in [S_F],\,k\in K,\,N_j\in\hat{S}_{F_1}\}$ for every node $N_j \in \hat{S}_F$;
\item either $X(\HH,L)$ is normal or else $S_F$ is of type $D_{m_0}+dA_h$ with  $n=m_0+d(h+1)$, $m_0\geq2$ and $h$ is even.
\end{enumerate}
\end{proposition}
\pf  Observe that (i) is the combinatorial translation of  \eqref{eq:condition} when $F'$ is a vertex of $\overline{F}$.  

If $\{kF~:~k\in K\}=\{F\}$, then $S_F$ cannot be of type  $D_{m_0}+dA_h$.  By Lemma \ref{lem:ferma} normality of  $X(\HH,L)$ is equivalent to normality of  $X$. If this is the case \eqref{eq:condition} holds whence (i) holds. 
We assume for the rest of the proof that $\#\{kF~:~k\in K\}>1$ and that $\#\hat{S}_F\geq2$, so $X$ is not a point. We set $B:=\hat{S}_F\cap \{N_0,N_1,N_{n-1},N_n\}$.

Assume first that $X$ is normal. By Lemma \ref{lem:idd}, condition (i) holds and  \eqref{eq:nodes} forces $\#B\neq3$. 

If $\#B=0$, then $\xi S_F=S_F$, so $\Sigma(\HH,L)\subseteq\Sigma(\HH,K,L)=\{F, \eta F\}$. Thus, either $\# \Sigma(\HH,L)=1$ and $X(\HH,L)$ is normal by Lemma \ref{lem:effe} and Theorem \ref{thm:lista-normali}, or 
$\# \Sigma(\HH,L)=2$. In this case Corollary \ref{cor:connected} implies that $\hat{S}_F$ and $\eta\hat{S}_F$ differ only by a vertex. Combined with (i) this is possible only if $\#\hat{S}_F=2$, $n=2m$ is even, and $N_m\in\hat{S}_F$. In this situation $X(\HH,L)$ is always normal, so  $\#\Sigma(\HH,L)=2$ cannot occur.

If $\#B=1$, we set $B=\{N_{\overline{j}}\}$. By assumption $X(\HH_{x_{\overline{j}}},K^{x_{\overline{j}}},L)=X(\HH_{x_{\overline{j}}},L)$ is normal and \cite[8.3.1]{richardson} forces $S_F$ to be is of type $D_{m_0}+dA_h$ with  $n=m_0+d(h+1)$, $m_0\geq2$, which gives (ii). 

If $\#B=4$, then $X(\HH_{x_j},K^{x_j},L)=X(\HH_{x_j},L)$ is normal for $j\in\{0,1,n-1,n\}$ but this contradicts \cite[8.3.1]{richardson} unless $S_F=\emptyset$. In this case $X(\HH,L)$ is normal.

Finally, assume $\# B=2$. If $B=\{N_0,N_1\}$, then $X(\HH_{x_0},K^{x_0},L)=X(\HH_{x_0},L)$ is normal, and  \cite[8.3.1]{richardson} forces $S_F$ to be of type $D_{m_0}$ so $X(\HH,L)$ is normal. 
 
The case $B=\{N_{n-1},N_n\}$ is dealt with in the same way. If $B=\{N_{\overline{i}},\,N_{\overline{j}}\}$ with $\overline{i}\leq1$ and $\overline{j}\geq n-1$, then $X(\HH_{x_{\overline{i}}},K^{x_{\overline{i}}},L)=X(\HH_{x_{\overline{i}}},L)$ and  \cite[8.3.1]{richardson} forces $S_F$ to be of type $dA_h$ with $h$ odd, so $X(\HH,L)$ is normal.

Conversely, assume that (i) and (ii) hold. In all cases $\#\Sigma(\HH,L)\leq 2$, so $X(\HH,L)$ is unibranch and $\Sigma(\HH,K,L)=\{kF~:~k\in K\}$. Condition (i) and Lemma \ref{lem:id} give thus  unibranchedness of $X$ at all vertices of faces in $\Sigma(\HH,K,L)$, i.e., at all minimal strata in $X$. It is also straightforward to verify that $X(\HH,L)$ is normal at all vertices with trivial stabiliser. By Corollary \ref{cor:soprasotto} the stratum $X$ is normal at such points. It remains to prove normality of $X(\HH_{x_j},K^{x_j},L)$ for every $j$ such that $x_j\in \overline{F}$ and $2\leq j\leq n-2$. The vertices of $k\overline{F}$ for any $k\in K$ are dealt with in a similar way.

If $n$ is odd or $n=2m$ and $j\neq m$, then $K^{x_j}=\langle \xi\rangle$, so $X(\HH_{x_j},K^{x_j},L)$ is also the finite counterpart at $x_j$ of the stratum corresponding to $S_F$ in $\mathrm{SO}_{2n}/\!/\mathrm{SO}_{2n}$,  which is normal by Proposition \ref{prop:so}. 

We finally look at $x_m$ for $n=2m$. Here, $K^{x_m}=K=\langle\eta\rangle\times\langle \xi\rangle$. 

If $S_F$ is of type $D_{m_0}\times dA_h\times D_{n_0}$ with $m_0\geq n_0\geq0$, $m_0\geq2$, then (i) applied to vertices not fixed by $\eta$ forces $m_0\geq m$. 
So if $N_{m}\in\hat{S}_F$ then $m=m_0$ and $X(\HH_{x_m},L)$ is a stratum for a finite Coxeter group of type $D_m\times D_m$ with trivial first component.  Hence,
\begin{align*}
X(\HH_{x_m},K,L)&= K\left(\{p\}\times X(\HH_{x_m,(2)},L_{(2)})\right)/K\\
&\simeq \langle\eta\rangle\left(\left(\{p\}\times \langle\xi\rangle X(\HH_{x_m,(2)},L_{(2)})\right)/\langle\xi\rangle\right)/\langle\eta\rangle\\
&\simeq \langle\eta\rangle\left(\{p\}\times X(\HH_{x_m},\langle\xi\rangle,L)\right)/\langle\eta\rangle\simeq \langle\eta\rangle\left(\{p\}\times X(\HH_B, L)\right)/\langle\eta\rangle.
\end{align*} 
The latter is normal by Example \ref{ex:dnbn}, from which we adopt notation, and  Example \ref{ex:giro}.

Finally, if $S_F$ is of type $dA_h$ with $h$ odd and $n=d(h+1)$, then, up to replacing $\xi$ by $\eta\xi$ and taking a graph automorphism of $\textrm{PSO}_{2n}$, we may always assume that $\eta S_F=S_F$, so $X(\HH_{x_m},K,L)\simeq X(\HH_{x_m},\langle\xi\rangle,L)/\langle\eta\rangle$ which is normal by Proposition \ref{prop:so}.
%\left(\{p\}\times X(\HH_{x_m,(2)},L_{(2)})\right)/K\\
%$X_K=\langle \xi,\eta \rangle X_\HH/\langle \xi,\eta \rangle\simeq\left(\langle \xi \rangle X_\HH)/\langle \xi \rangle\right)/\langle\eta\rangle$. The quotient $\langle \xi \rangle X_\HH/\langle \xi \rangle$ 
%is normal by Proposition \ref{prop:so}, so $X_K$ is normal.
\epf

\subsubsection{\bf  Type $E_6$ and $E_7$}

For type $E_6$ we only have the possibility $K=P^\vee/Q^\vee\simeq \mathbb{Z}/3\mathbb{Z}\simeq\langle \theta\rangle$ where the action of $\theta$ on the vertices of $\Abar$ is given by
$x_1\mapsto x_6\mapsto x_0\mapsto x_1$ and $x_3\mapsto x_5\mapsto  x_2\mapsto x_3$. The stratification in $\Ec/\WK$ corresponds to the Jordan stratification  in ${\mathrm E}_{6,{\mathrm {ad}}}/\!/{\mathrm E}_{6,{\mathrm{ad}}}$.

\begin{proposition}
Let $X=X(\HH,K,L)$ be a stratum in the Jordan stratification of ${\mathrm E}_{6,{\mathrm {ad}}}/\!/{\mathrm E}_{6,{\mathrm{ad}}}$ corresponding to the Coxeter class $[S_F]$.
Then $X$ is normal if and only if  $X(\HH,L)$ is normal.
%and \eqref{eq:nodes} holds.
%\comgio{vero? se funziona per $2A_2$ basta dire che e' normale se e solo se e' normale sopra}
\end{proposition}
\pf Assume first $\theta S_F=S_F$. If $X(\HH,L)$ is normal, then $X$ is normal by Lemma \ref{lem:ferma}. Conversely, if $X$ is normal,  then $X(\HH,L)$ is normal in codimension $1$, so by Proposition \ref{prop:codim1} and Theorem \ref{thm:lista-normali} either $X(\HH,L)$ is normal or  else $S_F$ is of type $D_4$. The latter is ruled out because
$(X(\HH,K,L),x_1)\sim_{se}(X(\HH_{x_1},K^{x_1},L),x_1)=(X(\HH_{x_1},L),x_1)$ which is not normal by \cite{broer,DR}. 
%In the latter case, $X_\HH$ is not normal at $x_1$. Corollary \ref{cor:soprasotto} at this vertex shows that $X_K$ is not normal. 
Also, the case of a point is immediate. We assume for the rest of the proof that $\theta S_F\neq  S_F$ and that $\#\hat{S}_F\geq2$. 

Assume $X$ is normal.  We verify by inspection that \eqref{eq:nodes} forces $S_F$ to be of type $A_5$, $2A_2$, $A_4+A_1$, or $A_2+2A_1$, $2A_2+A_1$. 

In the first two cases $X(\HH,L)$ is normal. The remaining cases cannot occur: indeed, the algorithm in Corollary \ref{cor:algo} shows that $\# \Sigma(\HH,K,L)\geq6$, contradicting Lemma \ref{lem:condition}.  
%so \eqref{eq:nodes} is not satisfied. 

Assume now that $X(\HH,L)$ is normal. Then, $S_F$ is either of type $A_5$ or $2A_2$. Observe that \eqref{eq:nodes} always holds. 
If $S_F$ is of type $A_5$, then $\dim X=1$ and it is normal in codimension $1$ by Proposition \ref{prop:equivalent}, hence it is normal. Let $S_F$ be of type $2A_2$. Lemma \ref{lem:id} gives unibranchedness of 
$X$ at all vertices of faces in $\Sigma(\HH,K,L)=\{F,\,\theta F,\,\theta^2F\}$. Corollary \ref{cor:soprasotto} implies that $X$ is normal at all $x_j\in \theta^i\overline{F}$ for $j\neq 4$ and $i\leq 3$. We finally apply Example \ref{ex:giro} to deduce normality of $X(\HH_{x_4},K^{x_4},L)$.
% using  Example \ref{ex:giro} because $X(\HH_{x_m},L)$  is the product of three  factors, cyclically permuted by $K$, and such that the only non-trivial one is normal.
\epf

\bigskip

For type $E_7$ we  have  $K=P^\vee/Q^\vee\simeq \mathbb{Z}/2\mathbb{Z}\simeq\langle k\rangle$ where $k$  fixes the vertices $x_2$ and $x_4$ and acts on the remaining ones by
$x_0\mapsto x_7\mapsto x_0$,  $x_1\mapsto x_6\mapsto x_1$, and $x_3\mapsto x_5\mapsto x_3$. The stratification in $\Ec/\WK$ corresponds to the Jordan stratification in ${\mathrm E}_{7,{\mathrm{ad}}}/\!/{\mathrm E}_{7,{\mathrm{ad}}}$.

\begin{proposition}
Let $X=X(\HH,K,L)$ be a stratum in the Jordan stratification of ${\mathrm E}_{7,{\mathrm{ad}}}/\!/{\mathrm E}_{7,{\mathrm{ad}}}$  corresponding to the Coxeter class $[S_F]$.
Then $X$ is normal if and only if $S_F$ satisfies one of the following conditions
\begin{enumerate}[label=(\roman*)]
\item it occurs in the list of Theorem \ref{thm:lista-normali} but it is not of type $D_4+A_1$ nor  $3A_1$;
\item it is of type $A_6$ or of type $A_5+ A_1$ and does not contain $N_2$. 
\end{enumerate}
\end{proposition}
\pf Observe that $S_F$ satisfies condition (i) if and only if  $X(\HH,L)$ is  normal and \eqref{eq:nodes} holds. 
Assume first $kS_F=S_F$ so (ii) does not occur. If $X(\HH,L)$ is normal, then $X$ is normal by Lemma \ref{lem:ferma} and \eqref{eq:nodes} trivially holds, so we have (i). Conversely, if $X$ is normal, then \eqref{eq:nodes} holds and Lemma \ref{lem:ferma} guarantees that $X(\HH,L)$ is normal in codimension $1$. Therefore, either  $S_F$ is as in (i) or it is of type $D_4$. This case can be excluded as we did for type $E_6$. Assume for the rest of the proof that $S_F\neq k S_F$ and that $\#\hat{S}_F\geq2$. 

If $X$ is normal,  by Lemma \ref{lem:allbutone} either $X(\HH,L)$ is normal in codimension $1$ or else $X(\HH,L)$ is unibranch and $\hat{S}_F$ contains at most three nodes, and at most one
different from $N_2$ and $N_4$. However, \eqref{eq:nodes} necessarily holds and under this assumption all strata $X(\HH,L)$ that are normal in codimension $1$ are normal, see  Proposition \ref{prop:codim1} and Theorem \ref{thm:lista-normali}, so we fall in (i). Assume that $\hat{S}_F$ contains at most three nodes, at most one different from $N_2$ and $N_4$ and that \eqref{eq:nodes} holds. By inspection, wee see that  $S_F$ is either as in (i) or (ii), or else  it is of type $A_3+A_2+A_1$, $A_2+A_4$, $A_3+A_2$ or $A_3+2A_1$.  Corollary \ref{cor:algo} applied to these $4$ cases would yield $\#\Sigma(\HH,K,L)>2$, contradicting Lemma \ref{lem:condition}.

Conversely, assume first  $S_F$ is as in (ii). Then $X$ is $1$-dimensional, with $\#\Sigma(\HH,K,L)\leq 2$ and satisfies the hypotheses of Proposition \ref{prop:equivalent}, so it is normal. The same argument gives normality for the strata $X$ for which $X(\HH,L)$ is $1$-dimensional and satisfies (i). There is only one case left satisfying (i) and $kS_F\neq S_F$, namely
$S_F$ of type $A_5$ containing $N_2$. Condition \eqref{eq:nodes} ensures unibranchedness of $X$ at all minimal strata. Since all nodes in $\hat{S}_F$ have trivial stabiliser, normality of $X$ follows from Corollary \ref{cor:soprasotto}. 
\epf

\section{Acknowledgements}
This work was partially supported by BIRD179758/17  and 
 DOR1691049/16 of the University of Padova. The first named author is a member of INdAM group GNSAGA.

\end{document}